\RequirePackage{fix-cm}
\documentclass[smallextended]{svjour3}       
\smartqed  
\usepackage{graphicx}
\usepackage{amsmath}	
\usepackage{amssymb}	
\usepackage[authoryear]{natbib}
\usepackage{hyperref}	
\hypersetup{colorlinks=true,linkcolor=blue,citecolor=blue,filecolor=blue,urlcolor=blue}
\newcommand{\mathbfit}[1]{\textbf{\textit{#1}}}
\newcommand{\mathbfss}[1]{\textbf{\textsf{#1}}}
\newcommand\degr{\hbox{$^\circ$}}
\newcommand{\doi}[1]{\textsc{doi}: \href{http://dx.doi.org/#1}{\nolinkurl{#1}}}
\journalname{Celestial Mechanics and Dynamical Astronomy}
\begin{document}

\title{Element sets for high-order Poincar\'e mapping of perturbed Keplerian motion}
%


\titlerunning{High-order Poincar\'e mapping of perturbed Keplerian motion}        

\author{David J. Gondelach \and Roberto Armellin}

\authorrunning{D. J. Gondelach, R. Armellin} 

\institute{D.J. Gondelach \at
              Surrey Space Centre, University of Surrey, GU2 7XH, Guildford, United Kingdom \\
              \email{d.gondelach@surrey.ac.uk}
           \and
           R. Armellin \at
              Surrey Space Centre, University of Surrey, GU2 7XH, Guildford, United Kingdom \\
              \email{r.armellin@surrey.ac.uk}
}

\date{Received: date / Accepted: date}

\maketitle

\begin{abstract}
The propagation and Poincar\'e mapping of perturbed Keplerian motion is a key topic in celestial mechanics and astrodynamics, e.g. to study the stability of orbits or design bounded relative trajectories. The high-order transfer map (HOTM) method enables efficient mapping of perturbed Keplerian orbits over many revolutions. For this, the method uses the high-order Taylor expansion of a Poincar\'e or stroboscopic map, which is accurate close to the expansion point. In this paper, we investigate the performance of the HOTM method using different element sets for building the high-order map. The element sets investigated are the classical orbital elements, modified equinoctial elements, Hill variables, cylindrical coordinates and Deprit's ideal elements. The performances of the different coordinate sets are tested by comparing the accuracy and efficiency of mapping low-Earth and highly-elliptical orbits perturbed by $J_2$ with numerical propagation. The accuracy of HOTM depends strongly on the choice of elements and type of orbit. A new set of elements is introduced that enables extremely accurate mapping of the state, even for high eccentricities and higher-order zonal perturbations. Finally, the high-order map is shown to be very useful for the determination and study of fixed points and centre manifolds of Poincar\'e maps.
\keywords{Poincar\'e map \and Stroboscopic map \and Orbit propagation \and Orbit element sets}
\end{abstract}


\section{Introduction}

The propagation of perturbed Keplerian motion is important for many different applications, such as 
predicting the orbit of a near-Earth satellite, studying the evolution of a planetary ring or designing a low-thrust trajectory. Propagating a perturbed orbit is, however, complicated by the fact that the differential equations of the dynamics are non-integrable and closed-form solutions cannot be obtained, except for special cases.
As a consequence, numerical integration is required to propagate the perturbed motion. 
Numerical propagation techniques, such as Cowell's method, are accurate, because they do not require any approximations of the dynamics \citep{vallado2013fundamentals}. However, they are also computationally inefficient, because accurate calculation of short-periodic effects requires the integration to be carried out using small time steps \citep{finkleman2014analysis}. Numerical propagation is therefore not suitable for long-term propagation or propagating many orbits, such as required for stability analysis.

The development of propagation methods that are both accurate and efficient is one of the key topics in celestial mechanics and astrodynamics.
In general, efficient propagation techniques rely on simplifying the dynamics. To this end, perturbation theory is applied by observing that perturbed Keplerian motion consists of pure Keplerian motion that is fast and integrable plus slow changes of the orbital plane and the in-plane motion caused by perturbations. These so-called general perturbation techniques can be divided in analytical and semi-analytical methods \citep{vallado2013fundamentals}.

Analytical methods use analytical approximations of the equations of motion, which allows analytical integration \citep{brouwer1959solution,kozai1962second}. These methods are extremely fast, however, due to the approximated dynamics they can only describe the characteristics of the motion for a limited time span \citep{vallado2013fundamentals}. Semi-analytical methods, on the other hand, combine numerical and analytical techniques 
to obtain a good combination of accuracy and efficiency. This is done by averaging the dynamics, i.e. filtering out the short-periodic motion, 
and integrating the {mean} orbital elements using large time steps \citep{deprit1969,mcclain1977recursively}. Deriving the averaged equations of motion and the equations for converting from osculating to mean elements can, however, be a complex task.

For the study of the evolution of a perturbed orbit, the continuous dynamics may not necessarily be of interest and a discrete dynamical system can be employed to compute the orbit at discrete moments in time. A well-known discrete method is the {Poincar\'e} map where the orbit state is mapped between two consecutive crossings of the orbit with a hyperplane, 
called a Poincar\'e section. Poincar\'e maps are frequently used to study the stability of quasi-periodic orbits and can elegantly show the orbital evolution of different orbits in the domain of interest, see e.g. \citet{Borderes2018}. An alternative to the Poincar\'e map is the {stroboscopic} map that maps the orbit over one period of the dynamics, e.g. at pericenter passage. 

Because the perturbed Keplerian dynamics are non-integrable a Poincar\'e or stroboscopic map has to be computed numerically unless simplifications are applied. 
Often semi-analytical techniques are employed to obtain Poincar\'e or stroboscopic maps.
For example, \citet{ely1997stability} apply averaging and Lie perturbation techniques to generate Poincar\'e plots for studying the stability of near-Earth orbits. \citet{Roth1978,roth1979} uses semi-analytical techniques to stroboscopically map a perturbed orbit from pericenter to pericenter to achieve efficient propagation. \citet{Broucke1994}, \citet{koon2001j2} and \citet{baresi2017}, on the other hand, rely on numerical propagation for mapping $J_2$-perturbed orbits to study their stability \citep{Broucke1994} and to find natural bounded relative trajectories \citep{koon2001j2,baresi2017}.

If the dynamics are approximately periodic, then consecutive points of Poincar\'e or stroboscopic mappings are close to each other. Consecutive evaluations of the map can therefore be approximated by a Taylor series expansion of the map at a preceding point.

\citet{berz1987} developed a technique to automatically compute high-order Taylor expansions of functions, called Taylor Differential Algebra (DA), and applied it to compute Taylor expansions of maps for particle beam dynamics. Recently, \citet{Wittig2015High} introducted this so-called high-order transfer map (HOTM) method in the field of astrodynamics and applied it to perturbed Keplerian motion of near-Earth satellites.
The HOTM is a high-order Taylor expansion of a Poincar\'e or stroboscopic map that is built by numerically propagating the orbit for one orbital revolution in Taylor differential algebra. This HOTM can be used to efficiently map the orbit over many revolutions. It was shown that the method allows one to accurately propagate orbits with reduced computation times compared to numerical propagation \citep{Wittig2015High,armellin2015hotm}.

In addition, because the map is built using numerical integration, any kind of perturbation can straightforwardly be included without the need for approximations.
Another advantage with respect to semi-analytical methods is that the osculating state is propagated and therefore no conversion from mean to osculating elements, and vice versa, is needed. As a result, the errors introduced by element conversions are omitted and the often complicated conversion equations do not have to be derived.

However, the main drawback of the HOTM propagation technique is the limited validity of the transfer map. The map consists of high-order Taylor expansions that are only accurate close the expansion point. Therefore, if the state that is mapped drifts away from the expansion point, the accuracy of the HOTM degrades. This characteristic of the HOTM also applies to changes in time in case of non-autonomous perturbations and therefore the HOTM accuracy reduces over time if the time dependency of the perturbations is not explicitly taken into account.

Generally, the accuracy and efficiency of propagation methods can be improved by selecting proper coordinate or orbital element sets for integration. 
This is achieved by rewriting the equations of motion in different variables in order to regularize and linearize the dynamics, in the sense of transforming nonlinear equations into linear ones without neglecting terms \citep{roa2017regularization}. Through the years, many different element sets have been developed and proposed for perturbed Keplerian motion \citep{hintz2008survey,bau2015elements}.
The HOTM method has only been implemented in classical orbital elements \citep{wittig2014long} and modified equinoctial elements \citep{Wittig2015High} without considering the impact of the choice of coordinates on the HOTM's validity.

In this paper, we investigate the HOTM performance for different coordinate sets. The element sets investigated are the classical orbital elements, modified equinoctial elements, Hill variables, cylindrical coordinates, Deprit's ideal elements and a novel set of elements. 
The performance of the HOTM for different coordinate sets is compared by propagating the orbit of a near-Earth satellite perturbed by the oblateness of the Earth. The causes for poor performance are analyzed and based on the results the new set of elements is introduced for improved performance. 
The best-performing coordinate set then is tested for higher-order zonal and drag perturbations. Finally, as an example application, the method is used to investigate quasi-periodic orbits around a fixed point.

Analysing the use of different element sets in this work has the goal of improving the HOTM method to extend its utility in orbital mechanics. This opens the door to new applications of the method to practical problems in astrodynamics or for theoretical studies in celestial mechanics. In addition, the method as presented in this paper is a starting point for extending it to different dynamical systems that involve other perturbations.

The paper is organised as follows. First, the applied dynamical model is discussed. After that, the different element sets and corresponding equations of motion are presented and their characteristics briefly compared. Then, the differential algebra technique and the high-order mapping method are introduced and the approach for building the high-order map is explained. The computation of fixed points of Poincar\'e maps is briefly discussed and the test cases are presented. Finally, the results are discussed and conclusions are drawn.

\section{Dynamical model}
The perturbations considered in this paper are zonal and drag perturbations. The fundamental equations for computing these perturbations are discussed in this section.

\subsection{Zonal perturbations}
\label{sec:geopotential}
The perturbing potential of an axially symmetric gravitational field is of the form \citep{wakker2015fundamentals}:
\begin{equation}
	R = \frac{\mu}{r} \sum^{\infty}_{n=2} J_n \left( \frac{R_e}{r} \right)^n P_n (\sin{\phi}),
\end{equation}
where $\mu$ is the gravitational parameter, $R_e$ is the equatorial radius, $J_n$ is the $n$-th zonal harmonic, $P_n (\sin{\phi})$ is the Legendre polynomial of degree $n$ in $\sin{\phi}$ and $\phi$ is the declination. The full potential is then given by:
\begin{equation}
	V = -\frac{\mu}{r} \left\{ 1-\sum^{\infty}_{n=2} J_n \left( \frac{R_e}{r} \right)^n P_n (\sin{\phi})  \right\} = -\frac{\mu}{r} + R. \label{eq:fullpotential}
\end{equation}
The acceleration due to the gravitational potential is obtained by taking the gradient of the potential function:
\begin{equation}
	\mathbfit{f} = -\boldsymbol{\nabla} V.
    \label{eq:gradient}
\end{equation}
When we consider only the second zonal harmonic, i.e. the $J_2$ term, the Legendre polynomial is:
\begin{equation}
	P_2(\sin{\phi}) = (3\sin^2{\phi}-1)/2,
    \label{eq:P2legendre}
\end{equation}
and the perturbing potential becomes:
\begin{equation}
	R_{J_2} = \frac{1}{2} \mu J_2 \frac{R_e^2}{r^3} (3 \sin^2{\phi} - 1).
    \label{eq:J2potential}
\end{equation}
The perturbing accelerations due to $J_2$ are then computed by taking the gradient: $\mathbfit{f}_{J_2} = -\boldsymbol{\nabla} R_{J_2}$.

\subsection{Drag}
\label{sec:drag}
The perturbing acceleration due to drag is given by \citep{vallado2013fundamentals}:
\begin{equation}
  \mathbfit{f}_{drag} = -\frac{1}{2} \rho\,C_d \frac{A}{m}\,|\mathbfit{V}_{rel}|\, \mathbfit{V}_{rel},
  \label{eq:dragPerturbation}
\end{equation}
where $C_d$ is the drag coefficient, $A/m$ the area-to-mass ratio, $\rho$ the atmospheric density and $\mathbfit{V}_{rel}$ the velocity vector with respect to the atmosphere. For density computations the axially-symmetric Harris-Priester atmospheric model \citep{harris1962time} is used.

\section{Element sets}
In this section, the element sets and corresponding equations of motion are introduced. In addition, the equations for computing the effect of the $J_2$ perturbation are provided. When possible we use equations of motion that are fully expressed in the elements used for the propagation to avoid introducing nonlinearities by converting between different coordinates. All elements are defined with respect to the Earth-centered inertial reference frame indicated by the axes $x$, $y$ and $z$.

\subsection{Classical orbital elements}
\begin{figure}
  \centering
	\includegraphics[width=0.7\columnwidth]{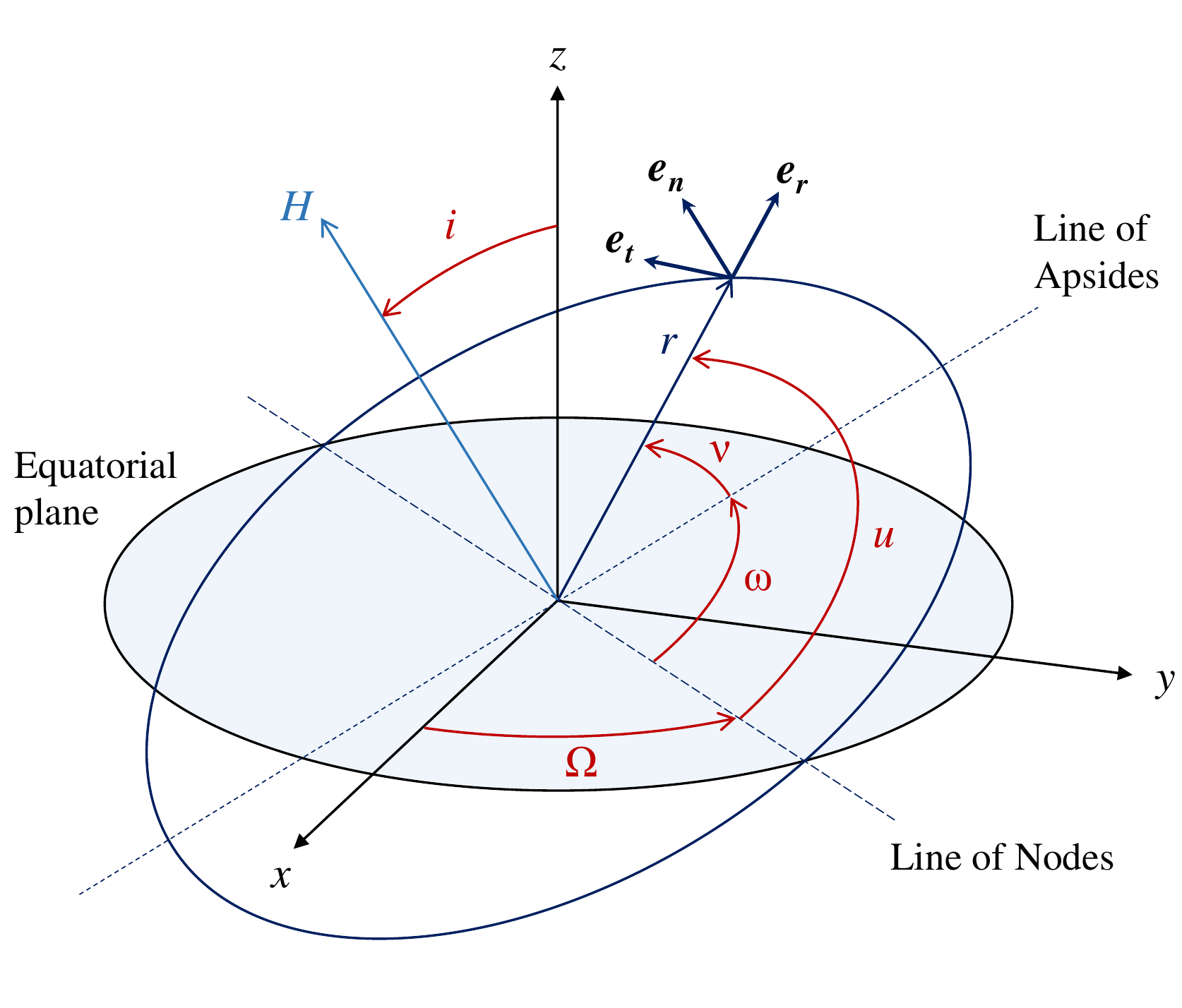}
    \caption{Diagram of orbital elements.}
    \label{fig:OrbitalElementsDrawing}
\end{figure}

The classical orbital elements (COE), as known as Keplerian elements, are given by \citep{vallado2013fundamentals}:
\begin{equation}
(a,e,i,\Omega,\omega,\nu),
\label{eq:Keplerianelements}
\end{equation}
where $a$ is the semi-major axis, $e$ the eccentricity, $i$ the inclination, $\Omega$ the right ascension of the ascending node, $\omega$ the argument of pericenter and $\nu$ the true anomaly, see Fig.~\ref{fig:OrbitalElementsDrawing}.

To compute the effect of a perturbing acceleration on the osculating Keplerian orbital elements one can apply Gauss' form of Lagrange's planetary equations \citep{battin1999introduction}:
\begin{align}
	\frac{da}{dt} &= \frac{2 a^2}{\sqrt{\mu p}} \left(e \sin{\nu} f_r +\frac{p}{r} f_t \right) ,
	\label{eq:dadt} \\
	\frac{de}{dt} &= \frac{1}{\sqrt{\mu p}} \left[ p \sin{\nu} f_r + \{(p+r)\cos{\nu} + e\,r \} f_t \right] ,
	\label{eq:dedt} \\
	\frac{di}{dt} &= \frac{r \cos{u}}{\sqrt{\mu p}} f_n ,
	\label{eq:didt} \\
	\frac{d\Omega}{dt} &= \frac{r \sin{u}}{\sqrt{\mu p} \sin{i}} f_n ,
	\label{eq:dRAANdt} \\
	\frac{d\omega}{dt} &= \frac{1}{e \sqrt{\mu p}} \left[ - p \cos{\nu} f_r + (p+r) \sin{\nu} f_t \right]-\frac{r \cos{i} \sin{u}}{\sqrt{\mu p} \sin{i}} f_n ,
	\label{eq:dargPerdt} \\
    \frac{d\nu}{dt} &= \frac{\sqrt{\mu p}}{r^2} + \frac{1}{e\sqrt{\mu p}} \left[ p \cos{\nu} f_r - (p+r) \sin{\nu} f_t \right] ,
\end{align} \\
where $f_r$, $f_{t}$ and $f_n$ are the components of the perturbing acceleration in the radial, transverse and normal directions\footnote{The radial direction $\mathbfit{e}_r$ points along the radius vector, the transverse direction $\mathbfit{e}_t$ is normal to the radial direction in the orbital plane and the normal direction $\mathbfit{e}_n$ is normal to the orbital plane along the angular momentum vector, see Fig.~\ref{fig:OrbitalElementsDrawing}.}, respectively. 
In addition, $u$ is the argument of latitude:
\begin{equation}
u=\omega + \nu ,
\label{eq:argumentOfLatitude}
\end{equation}
$p$ is the semi-latus rectum:
\begin{equation}
p = a(1-e^2) ,
\label{eq:semiLatusRectum}
\end{equation}
and $r$ is the radial distance:
\begin{equation}
r = \frac{p}{1+e\cos{\nu}} .
\end{equation}

The perturbing forces can be computed using the perturbing potential expressed in spherical coordinates $R=R(r,u,i)$ and taking the gradient \citep{vallado2013fundamentals}:

\begin{equation}
f_r = -\frac{\partial R}{\partial r}, ~~~~~~~~ f_t = -\frac{1}{r} \frac{\partial R}{\partial u}, ~~~~~~~~ f_n = -\frac{1}{r\sin{u}} \frac{\partial R}{\partial i} .
\label{eq:RTNgradient}
\end{equation} \\
Considering that $\sin{\phi} = \sin{i} \sin{u}$, we obtain the $J_2$ perturbing forces as:
\begin{align}
	f_r &= \frac{3}{2} \mu J_2 \frac{R_e^2}{r^4} \left( 3 \sin^2{i} \sin^2{u}-1 \right) ,
	\label{eq:J2radial} \\
	f_t &= -3 \mu J_2 \frac{R_e^2}{r^4} \sin^2{i} \sin{u} \cos{u} ,
	\label{eq:J2transverse} \\
	f_n &= -3 \mu J_2 \frac{R_e^2}{r^4} \sin{i} \cos{i} \sin{u} .
	\label{eq:J2normal}
\end{align}

\subsection{Modified equinoctial elements}
The modified equinoctial elements (MEE) are defined by \citet{Walker1985} as:
\begin{equation}
\begin{array}{lll}
p = a(1-e^2),					&	f = e\cos{(\omega+\Omega)},		&	g = e\sin{(\omega+\Omega)}, \\
h = \tan{(i/2)}\cos{\Omega},~	&	k = \tan{(i/2)}\sin{\Omega},~	&	L = \Omega+\omega+\nu .
\end{array} 
\label{eq:equinoctialelements}
\end{equation}
Lagrange's planetary equations can be written for the modified equinoctial elements as follows \citep{Walker1985,Walker1986}:

\begin{align}
\frac{dp}{dt} =& 2 \sqrt{\frac{p}{\mu}} \left(-g \frac{\partial R}{\partial f} + f \frac{\partial R}{\partial g} + \frac{\partial R}{\partial L} \right) ,
\label{eq:dpdt}	\\
	\frac{df}{dt} =& \frac{1}{\sqrt{\mu p}} \left \{ 2pg\frac{\partial R}{\partial p} - \left(1-f^2-g^2\right)\frac{\partial R}{\partial g} - \frac{g s^2}{2} \left(h \frac{\partial R}{\partial h}+k\frac{\partial R}{\partial k} \right) \right.
	\notag	\\
	&+ \left. \left[f+(1+w)\cos{L} \right]\frac{\partial R}{\partial L} \right \} ,
	\label{eq:dfdt}	\\
	\frac{dg}{dt} =& \frac{1}{\sqrt{\mu p}} \left \{ -2pf\frac{\partial R}{\partial p} + \left(1-f^2-g^2\right)\frac{\partial R}{\partial f} + \frac{f s^2}{2} \left(h\frac{\partial R}{\partial h}+k\frac{\partial R}{\partial k}\right) \right.
	\notag	\\
	&+ \left. \left[g+(1+w)\sin{L}\right]\frac{\partial R}{\partial L} \right \} ,
	\label{eq:dgdt}	\\
	\frac{dh}{dt} =& \frac{s^2}{2\sqrt{\mu p}} \left \{h \left(g\frac{\partial R}{\partial f} - f\frac{\partial R}{\partial g} - \frac{\partial R}{\partial L} \right) - \frac{s^2}{2}\frac{\partial R}{\partial k} \right \} ,
	\label{eq:dhdt}	\\
	\frac{dk}{dt} =& \frac{s^2}{2\sqrt{\mu p}} \left \{k \left(g\frac{\partial R}{\partial f} - f\frac{\partial R}{\partial g} - \frac{\partial R}{\partial L} \right) + \frac{s^2}{2}\frac{\partial R}{\partial h} \right \} ,
	\label{eq:dkdt}	\\
	\frac{dL}{dt} =& \sqrt{\mu p}\left( \frac{w}{p} \right)^2 + \frac{s^2}{2\sqrt{\mu p}} \left\{ h\frac{\partial R}{\partial h} + k\frac{\partial R}{\partial k} \right \} ,
	\label{eq:dLdt}
\end{align}
where \( w = 1+f\cos{L}+g\sin{L} \) and \( s^2 = 1+h^2+k^2 \).

By writing the $J_2$ perturbing potential \eqref{eq:J2potential} in modified equinoctial elements using $r=p/w$ and
\begin{equation}
	\sin{\phi} = \frac{2(h\sin{L}-k\cos{L})}{s^2} ,
	\label{eq:sinphiMEE}
\end{equation}
the partial derivatives of $R$ can be computed as:
\begin{align}
	\frac{\partial R}{\partial p} =& \frac{3\mu}{w r^2}~ J_2 \left(\frac{R_e}{r}\right)^2 P_2(\sin{\phi}) ,
	\label{eq:dRdp}	\\
	\frac{\partial R}{\partial f} =& \frac{-3\mu \cos{L}}{w r}~ J_2 \left(\frac{R_e}{r}\right)^2 P_2(\sin{\phi}) ,
	\label{eq:dRdf}	\\
	\frac{\partial R}{\partial g} =& \frac{-3\mu \sin{L}}{w r}~ J_2 \left(\frac{R_e}{r}\right)^2 P_2(\sin{\phi}) ,
	\label{eq:dRdg}	\\
	\frac{\partial R}{\partial h} =& \frac{-2 \mu}{r s^4}   \left\{ (1-h^2+k^2) \sin{L} + 2 h k \cos{L}\right\} ~ J_2 \left(\frac{R_e}{r}\right)^2   P_2'(\sin{\phi}) ,
	\label{eq:dRdh}	\\
	\frac{\partial R}{\partial k} =& \frac{2 \mu}{r s^4}   \left\{(1+h^2-k^2) \cos{L} + 2 h k \sin{L} \right\} ~ J_2 \left(\frac{R_e}{r}\right)^2   P_2'(\sin{\phi}) ,
	\label{eq:dRdk}	\\
	\frac{\partial R}{\partial L} =& \frac{-2 \mu}{r s^2}   (h \cos{L} + k \sin{L}) ~ J_2 \left(\frac{R_e}{r}\right)^2   P_2'(\sin{\phi})   ,
	\notag  \\
	& - \frac{3\mu}{w r}   (g \cos{L}-f \sin{L})~ J_2 \left(\frac{R_e}{r}\right)^2   P_2(\sin{\phi}) ,
	\label{eq:dRdL}
\end{align}\\
with $P_2$ from Eq.~\eqref{eq:P2legendre} and $P_2'(\sin{\phi}) = \frac{dP_2}{d(\sin{\phi})} = 3\sin{\phi}$.

The Gauss' equations of motion in terms of modified equinoctial elements are given in Appendix~\ref{app:EoMMEE}.

\subsection{Hill variables}
The Hill variables, also known as polar-nodal variables or Whittaker variables, are canonical variables and defined as\footnote{The Hill variables are often written as $(r, \theta, \nu, R, \Theta, N)$ where $\theta = u$, $\nu = \Omega$, $R = \dot{r}$, $\Theta=H$ and $N=H_z$.} \citep{hill1913}:
\begin{equation}
(r, u, \Omega, \dot{r}, H, H_z) ,
\label{eq:polarnodalelements}
\end{equation}
where $r$ is the radial distance, $u$ is the argument of latitude, $\Omega$ is the right ascension
of the ascending node, $\dot{r} = dr/dt$ is the radial velocity, $H$ is the angular momentum, and $H_z = H \cos{i}$ is the polar component of the angular momentum.

Because the Hill variables are canonical variables, the equations of motion can be obtained directly from the Hamiltonian. The Hamiltonian for an autonomous conservative system is the sum of the potential and kinetic energy. Using the potential \eqref{eq:fullpotential} where $\sin{\phi} = \sin{i} \sin{u}$ and $\sin^2{i} = 1-H_z^2/H^2$, we can write the Hamiltonian as follows:
\begin{equation}
	\mathcal{H} = \frac{1}{2} \left({\dot{r}^2} + \frac{{H}^2}{r^2}\right) - \frac{\mu}{r} \left\{1 - \frac{1}{2} J_2 \frac{{R_e}^2}{r^2} \left[3\, {\sin^2u}\, \left( 1 - \frac{H_z^2}{{H}^2} \right) - 1\right] \right\} ,
\end{equation}
where $H/r$ is the transverse velocity.

The equations of motion (i.e. Hamilton's equations) are then obtained as:
\begin{align}
	\frac{dr}{dt} &= \frac{\partial \mathcal{H}}{\partial \dot{r}} = \dot{r} ,
	\label{eq:drdtpolarnodal} 	\\
	\frac{du}{dt} &= \frac{\partial \mathcal{H}}{\partial H} = \frac{H}{r^2} + \frac{3\, \mu\, J_2\, {R_e}^2\, H_z^2\, {\sin^2u}}{r^3 {H}^3} ,
	\label{eq:dudtpolarnodal}	\\
	\frac{d\Omega}{dt} &= \frac{\partial \mathcal{H}}{\partial H_z} = -\frac{3\, \mu\, J_2\, {R_e}^2\, H_z\, {\sin^2u}}{r^3 {H}^2} ,
	\label{eq:dRAANdtpolarnodal}	\\
	\frac{d\dot{r}}{dt} &= -\frac{\partial \mathcal{H}}{\partial r} = 
	\frac{H^2}{r^3} -\frac{\mu }{r^2} + \frac{3}{2} \mu J_2 \frac{R_e^2}{r^4} \left[ 3\sin^2{u} \left(1-\frac{H_z^2}{H^2}\right) -1 \right] ,
	\label{eq:drdotdtpolarnodal}	\\
	\frac{dH}{dt} &= -\frac{\partial \mathcal{H}}{\partial u} = -3\, \mu\, J_2\, \frac{R_e^2}{r^3}\, \cos{u}\, \sin{u}\, \left(1 - \frac{H_z^2}{{H}^2} \right) ,
	\label{eq:dHdtpolarnodal}	\\
	\frac{dH_z}{dt} &= -\frac{\partial \mathcal{H}}{\partial \Omega} = 0 .
	\label{eq:dHzdtpolarnodal}
\end{align}
For completeness, the Gaussian form of the equations of motion in Hill variables is given in Appendix~\ref{app:EoMpolarnodal}.

\subsection{Cylindrical coordinates}
\label{sec:cylcoord}

\begin{figure}
  \centering
	\includegraphics[width=0.7\columnwidth]{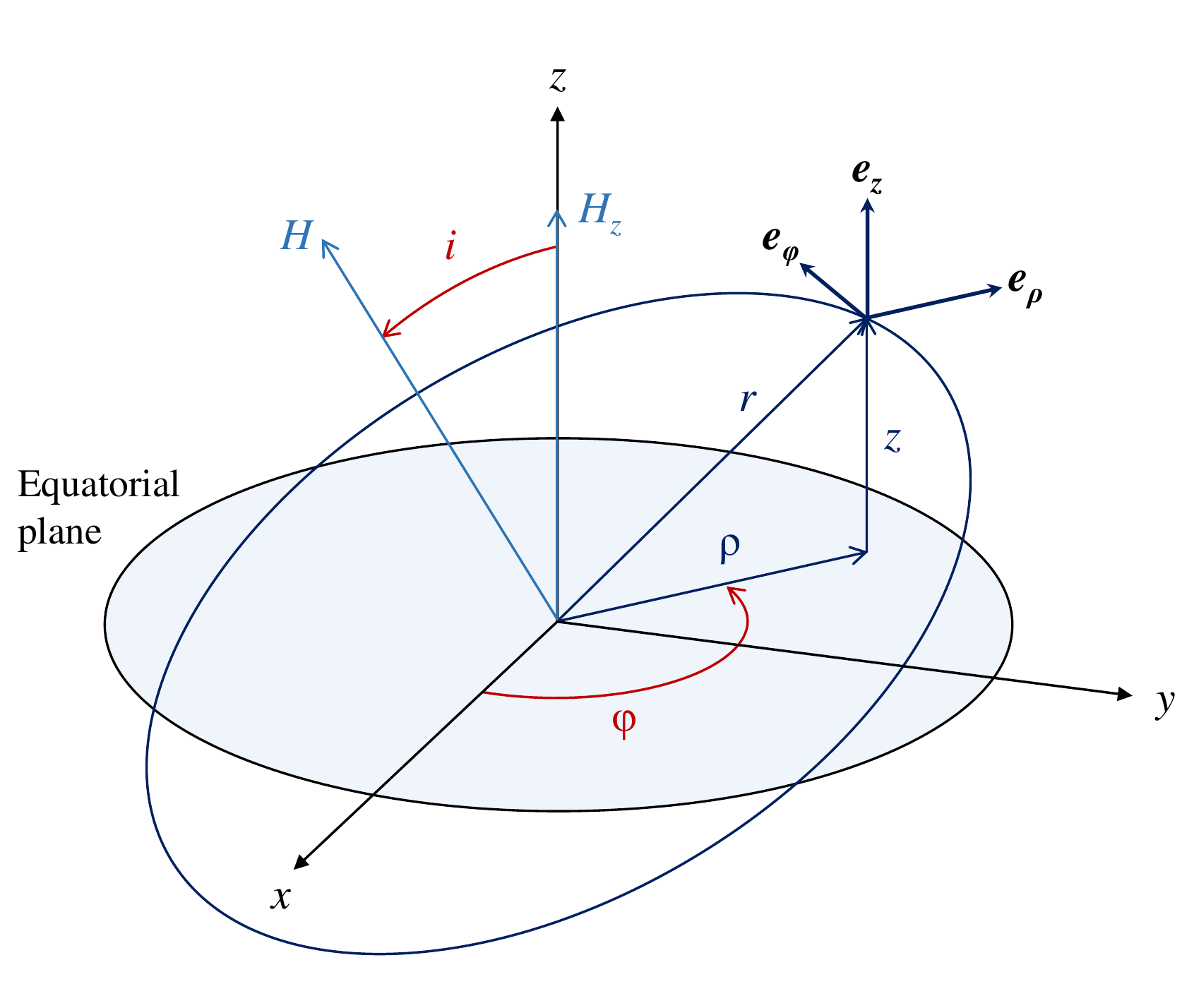}
    \caption{Diagram of cylindrical coordinates.}
    \label{fig:CylindricalCoordinatesDrawing}
\end{figure}

The cylindrical coordinates $(\rho,\varphi,z)$ are defined by the distance from the $z$-axis $\rho$, the azimuth angle $\varphi$ and the height $z$, see Fig.~\ref{fig:CylindricalCoordinatesDrawing}. 
Together with their time derivatives $(\dot{\rho},\dot{\varphi},\dot{z})$ the coordinate set (Cyl) can be used for orbit propagation. 
The equations of motion in terms of cylindrical coordinates are:
\begin{align}
	\frac{d\rho}{dt} &= \dot{\rho} , \label{eq:drhodt} \\
	\frac{d\varphi}{dt} &= \dot{\varphi} , \label{eq:dvarphidt} \\
	\frac{dz}{dt} &= \dot{z} , \\
	\frac{d\dot{\rho}}{dt} &= \rho\dot{\varphi}^2 + f_{\rho} , \label{eq:drhodotdt} \\
	\frac{d\dot{\varphi}}{dt} &= \frac{-2\dot{\rho}\dot{\varphi}}{\rho} + \frac{1}{\rho} f_{\varphi} , \label{eq:dvarphidotdt} \\
	\frac{d\dot{z}}{dt} &= f_z , \label{eq:ddotzdt}
\end{align}
where $f_{\rho}$, $f_{\varphi}$ and $f_{z}$ are the forces in projected radial, azimuthal and axial direction. These forces are computed by taking the gradient of the full potential \eqref{eq:fullpotential} expressed in cylindrical coordinates using \( r =\sqrt{\rho^2+z^2} \) and $\sin{\phi} = z/r$:
\begin{align}
	f_{\rho} &= -\frac{\partial V}{\partial \rho} = -\frac{\mu}{r^3}\rho + \frac{1}{2}\frac{\rho}{r^7} J_2 \mu R_e^2 (12 z^2 - 3 \rho^2) ,
	\label{eq:fJ2rho} \\
	f_{\varphi} &= -\frac{1}{\rho}\frac{\partial V}{\partial \varphi} = 0 ,
	\label{eq:fJ2varphi} \\
	f_{z} &= -\frac{\partial V}{\partial z} = -\frac{\mu}{r^3}z + \frac{1}{2}\frac{z}{r^7} J_2 \mu R_e^2 (6 z^2 - 9 \rho^2) .
	\label{eq:fJ2z}
\end{align}
These forces include the Keplerian part of the gravitational attraction that is not incorporated in equations of motion \eqref{eq:drhodt}-\eqref{eq:ddotzdt}.

As an alternative to the angular velocity $\dot{\varphi}$, the z-component of the angular momentum $H_z$ can be used and the element set (CylHz) becomes $(\rho,\varphi,z,\dot{\rho},\dot{z},H_z)$. The equations of motion involving $\dot{\varphi}$ (\eqref{eq:dvarphidt}, \eqref{eq:drhodotdt} and \eqref{eq:dvarphidotdt}) are then replaced by \citep{Deprit1994Linearization}:
\begin{align}
	\frac{d\varphi}{dt} &= \frac{H_z}{\rho^2} , \\
	\frac{d\dot{\rho}}{dt} &= \frac{H_z^2}{\rho^3} + f_{\rho} , \\
	\frac{dH_z}{dt} &= \rho f_{\varphi} .
\end{align}

\subsection{Ideal elements}
\begin{figure}
  \centering
	\includegraphics[width=0.7\columnwidth]{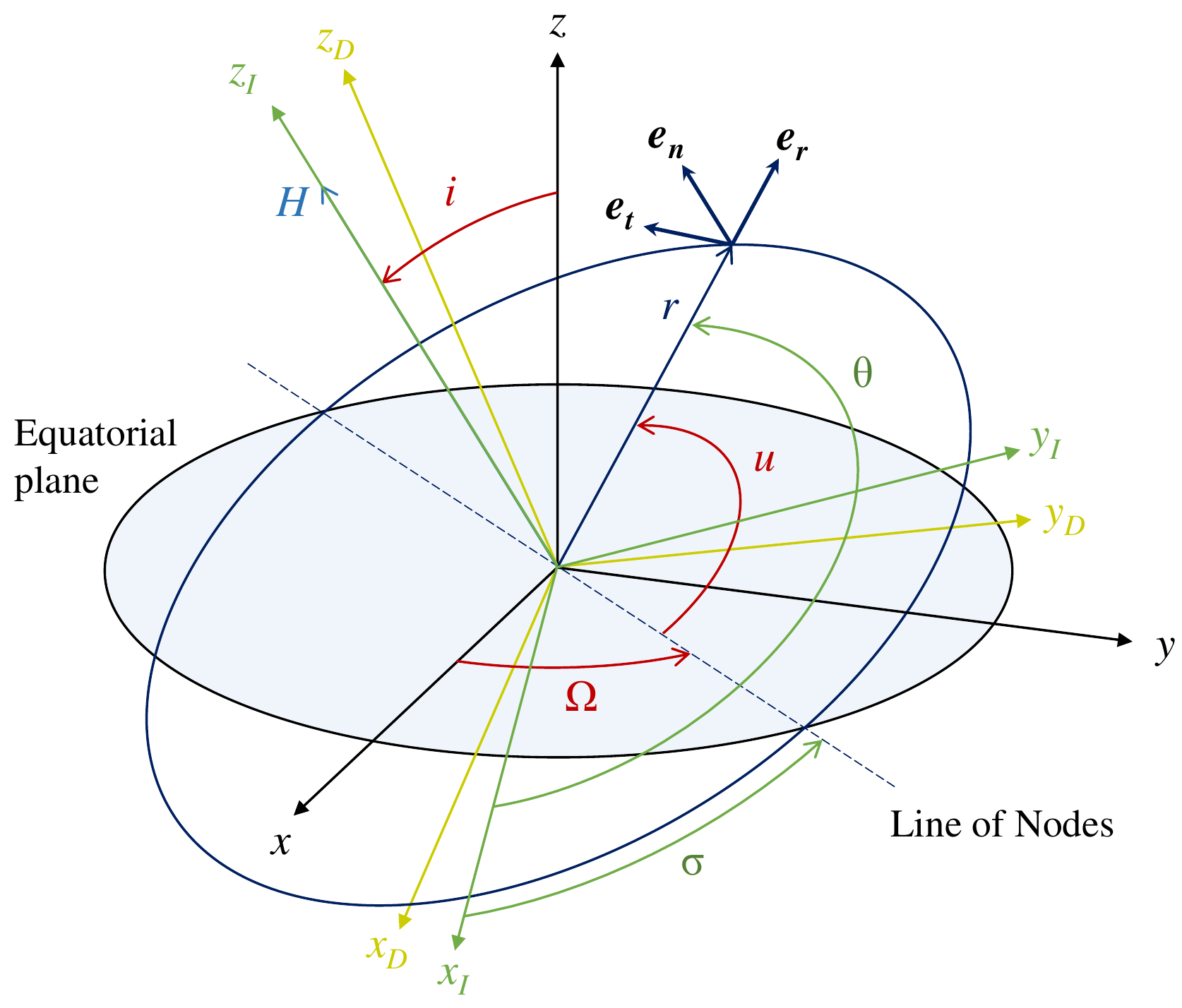}
    \caption{Diagram of the ideal elements and the ideal frame $\mathcal{I}$ and departure frame $\mathcal{D}$. The orientation of $\mathcal{D}$ with respect to the inertial frame is defined by the rotation matrix $\mathbfss{M}$ and the orientation of $\mathcal{I}$ with respect to $\mathcal{D}$ by the quaternion $\boldsymbol{\lambda}$.}
    \label{fig:IdealElementsDrawing}
\end{figure}

The ideal elements were developed by \citet{deprit1975ideal} and use the concept of Hansen's ideal frame \citep{hansen1857auseinandersetzung} to (partially) decouple the fast in-plane motion from the slow rotation of the orbital plane. The orientation of the orbital plane is determined by the ideal frame $\mathcal{I}$ that rotates slowly with respect to the departure frame $\mathcal{D}$ and is defined by a quaternion $\boldsymbol{\lambda}$, see Fig.~\ref{fig:IdealElementsDrawing}.
The attitude of the departure frame in the inertial frame is given by the initial right ascension of the ascending node $\Omega_0$, inclination $i_0$ and argument of latitude $u_0$, that define the rotation matrix $\mathbfss{M}$.

The ideal elements developed by \citet{deprit1975ideal} are defined as:
\begin{equation}
(\lambda_1, \lambda_2, \lambda_3, \lambda_4, H, C, S, \theta) ,
\label{eq:idealelements}
\end{equation}
where $\theta$ is the in-plane angle of the position vector in the ideal frame, $H$ is the angular momentum, and $C$ and $S$ are related to the direction of the eccentricity vector in the ideal frame and are given by:
\begin{align}
	C &= \left( \frac{H}{r} - \frac{H}{p} \right) \cos{\theta} + \dot{r} \sin{\theta} ,
	\label{eq:IdealC} \\
	S &= \left( \frac{H}{r} - \frac{H}{p} \right) \sin{\theta} - \dot{r} \cos{\theta} .
	\label{eq:IdealS}
\end{align}
Finally, $\lambda_i$ are the components of the quaternion that relate the ideal frame to the departure frame:
\begin{align}
	\lambda_1 &= \sin{ \left(\tfrac{1}{2}i_{\mathcal{I}} \right)} \cos{\left( \tfrac{1}{2}(\Omega_{\mathcal{I}} - \sigma_{\mathcal{I}}) \right)} ,
	\label{eq:IdealLambda1} \\
	\lambda_2 &= \sin{ \left(\tfrac{1}{2}i_{\mathcal{I}} \right)} \sin{\left( \tfrac{1}{2}(\Omega_{\mathcal{I}} - \sigma_{\mathcal{I}}) \right)} ,
	\label{eq:IdealLambda2} \\
	\lambda_3 &= \cos{ \left(\tfrac{1}{2}i_{\mathcal{I}} \right)} \sin{\left( \tfrac{1}{2}(\Omega_{\mathcal{I}} + \sigma_{\mathcal{I}}) \right)} ,
	\label{eq:IdealLambda3} \\
	\lambda_4 &= \cos{ \left(\tfrac{1}{2}i_{\mathcal{I}} \right)} \cos{\left( \tfrac{1}{2}(\Omega_{\mathcal{I}} + \sigma_{\mathcal{I}}) \right)} ,
	\label{eq:IdealLambda4}
\end{align}
where $\Omega_{\mathcal{I}}$, $\sigma_{\mathcal{I}}$ and $i_{\mathcal{I}}$ denote the corresponding Euler angles.

%
%

The equations of motion of the ideal elements are given by \citep{lara2017ideal}:
\begin{align}
	\dot{\lambda_1} &= \frac{r}{2H} f_n (\lambda_4\cos{\theta} - \lambda_3\sin{\theta}) ,
	\label{eq:IdealdLambda1dt} \\
	\dot{\lambda_2} &= \frac{r}{2H} f_n (\lambda_4\sin{\theta} + \lambda_3\cos{\theta}) ,
	\label{eq:IdealdLambda2dt} \\
	\dot{\lambda_3} &= \frac{r}{2H} f_n (\lambda_1\sin{\theta} - \lambda_2\cos{\theta}) ,
	\label{eq:IdealdLambda3dt} \\
	\dot{\lambda_4} &= \frac{r}{2H} f_n (-\lambda_1\cos{\theta} - \lambda_2\sin{\theta}) ,
	\label{eq:IdealdLambda4dt} \\
	\dot{H} &= r\,f_t ,
	\label{eq:IdealdHdt} \\
	\dot{C} &= \left( 1 + \frac{r}{p} \right) f_t \cos{\theta} + f_r \sin{\theta} ,
	\label{eq:IdealdCdt} \\
	\dot{S} &= \left( 1 + \frac{r}{p} \right) f_t \sin{\theta} - f_r \cos{\theta} ,
	\label{eq:IdealdSdt} \\
	\dot{\theta} &= \frac{H}{r^2} ,
	\label{eq:Idealdthetadt}
\end{align} \\
where $p=H^2/\mu$ and 
\begin{equation}
	r = \frac{H}{\tfrac{H}{p}+C\cos{\theta}+S\sin{\theta}} .
\end{equation}
The perturbing forces in the orbital frame $(f_r,f_t,f_n)$ are computed by calculating the perturbations in the inertial frame and transforming them to the orbital frame via the departure and ideal frames using three rotations defined by the rotation matrix $\mathbfss{M}$, quaternion $\boldsymbol{\lambda}$ and ideal angle $\theta$ \citep{lara2017ideal}.

The $J_2$ perturbations in the inertial $x$, $y$ and $z$ directions are obtained by writing the $J_2$ potential \eqref{eq:J2potential} in Cartesian coordinates using $r=\sqrt{x^2+y^2+z^2}$ and $\sin{\phi}=z/r$ and taking the gradient:
\begin{align}
	f_x &= -\frac{\partial R}{\partial x} = -\frac{3}{2} \mu J_2 \frac{R_e^2}{r^5} x \left( 1 - 5 \frac{z^2}{r^2} \right) , \\
	f_y &= -\frac{\partial R}{\partial y} = -\frac{3}{2} \mu J_2 \frac{R_e^2}{r^5} y \left( 1 - 5 \frac{z^2}{r^2} \right) , \\
	f_z &= -\frac{\partial R}{\partial z} = -\frac{3}{2} \mu J_2 \frac{R_e^2}{r^5} z \left( 3 - 5 \frac{z^2}{r^2} \right) .
\end{align}

\subsection{Eccentric Hill variables}
\label{sec:axialnodal}
Finally, a new set of orbital elements is introducted to improve to performance of the HOTM. This new set is called the eccentric Hill variables (EccHill), because they are closely related to the Hill variables, but use eccentric variables instead of $r$ and $\dot{r}$. The eccentric Hill variables are defined as\footnote{It can be noticed that this set of parameters is very similar to the elements used by \citet{deprit1970} to obtain an analytical solution for the main problem in satellite theory: $(\sqrt{\mu a},H_z,\hat{f},\hat{g},\Omega,M+\omega)$.}: 
\begin{equation}
(H,H_z,\hat{f},\hat{g},\Omega,u) ,
\label{eq:axialnodalelements}
\end{equation}
where $\hat{f}$ and $\hat{g}$ are components of the eccentricity vector:
\begin{align}
	\hat{f} &= e\cos{\omega} ,	\\
	\hat{g} &= e\sin{\omega} .
\end{align}
The time derivatives of $\hat{f}$ and $\hat{g}$ are obtained by observing that:
\begin{align}
	\frac{d\hat{f}}{dt} &= \frac{de}{dt} \cos{\omega} - e \sin{\omega} \frac{d\omega}{dt} ,	\\
	\frac{d\hat{g}}{dt} &= \frac{de}{dt} \sin{\omega} + e \cos{\omega} \frac{d\omega}{dt} ,
\end{align}
and taking $de/dt$ and $d\omega/dt$ from equations \eqref{eq:dedt} and \eqref{eq:dargPerdt}.

The equations of motion in terms of eccentric Hill variables are then obtained as:
\begin{align}
	\frac{dH}{dt} &= r f_t ,	\\
	\frac{dH_z}{dt} &= \frac{r}{H} \left( H_z f_t - G\cos{u} f_n \right) ,	\\
	\frac{d\hat{f}}{dt} &= \frac{r}{H} \left\{ \hat{w}\sin{u} f_r + \left[(\hat{w}+1) \cos{u}+\,\hat{f}\right] f_t + \frac{\hat{g}\,H_z\sin{u}}{G} f_n\right\} ,	\\
	\frac{d\hat{g}}{dt} &= \frac{r}{H} \left\{-\hat{w}\cos{u} f_r + \left[(\hat{w}+1) \sin{u}+\hat{g}\right] f_t - \frac{\hat{f}\,H_z\sin{u}}{G} f_n\right\} , 	\\
	\frac{d\Omega}{dt} &= \frac{r\sin{u}}{G} f_n ,	\\
	\frac{du}{dt} &= \frac{H}{r^2} - \frac{r H_z \sin{u}}{H G} f_n ,
\end{align}
where $G=\sqrt{H^2-H_z^2}$, $\hat{w}=1+\hat{f}\cos{u}+\hat{g}\sin{u}$ and $r=H^2/(\mu \hat{w})$.

The $J_2$ perturbing forces can be obtained directly from equations \eqref{eq:J2radial}-\eqref{eq:J2normal} using $\cos{i}=H_z/H$ and $\sin{i} = G/H$. Alternatively, the potential $R_{J_2}$ can be expressed in eccentric Hill variables to compute $f_n$ as:
\begin{equation}
f_n = -\frac{1}{r\sin{u}} \frac{\partial R}{\partial H_z} \frac{\partial H_z}{\partial i} = \frac{H \sin{i}}{r\sin{u}} \frac{\partial R}{\partial H_z} = \frac{G}{r\sin{u}} \frac{\partial R}{\partial H_z} .
\label{eq:fhAxialNodal}
\end{equation}

\subsection{Comparison}
The intrinsic properties of the various element sets make different coordinates more suitable for propagation than others depending on the orbit and the dynamics. Particularly, the use of specific element sets may result in singularities in the equations of motion. If the eccentricity is zero, the argument of perigee $\omega$ is not defined and the COE become singular. Similarly, when the inclination is zero, the longitude of the node $\Omega$ is undefined and the COE, Hill and EccHill variables become singular. The cylindrical coordinates, on the other hand, contain a singularity at the poles, i.e. when the inclination is 90$\degr$, where the azimuth angle $\varphi$ is not defined and the distance $\rho$ is zero. The CylHz coordinates may, however, be used at $i=90\degr$ because $\rho$ only vanishes when $H_z$ also vanishes, thus canceling the singularity in the equations of motion. None of the coordinate sets contains a singularity at the critical inclination, $i=63.4\degr$, which sometimes causes singularities in analytical solutions to perturbed Keplerian motion \citep{brouwer1959solution}. An overview of the singularities is shown in Table~\ref{tab:CoordinateComparison}. 

Besides, for zonal perturbations the dynamics do not depend on the longitude, $\Omega$ and $\varphi$ (see e.g. \eqref{eq:dHzdtpolarnodal} and \eqref{eq:fJ2varphi}), and consequently the polar component of the angular momentum $H_z$ is constant.

\begin{table*}
  \centering
  \caption{Overview of the element sets, their singularities and the element used as independent variable for propagation.}
    \begin{tabular}{lccc}
    \hline
    \textbf{Element set} & \textbf{Elements} & \textbf{Independent} & \multicolumn{1}{c}{\textbf{Singularities}} \\
    \textbf{} & & \textbf{variable} &  \\
    \hline
    \textbf{COE} 	& $a,e,i,\Omega,\omega,\nu$ & $\nu$ 	& ${e=0}$, ${i=0\degr,  180\degr}$ \\
    \textbf{MEE} 	& $p,f,g,h,k,L$ & $L$     	& $-$     \\
    \textbf{Cyl} 	& $\rho,\varphi,z,\dot{\rho},\dot{\varphi},\dot{z}$ & $t$     	& ${i=90\degr}$  \\
    \textbf{CylHz} 	& $\rho,\varphi,z,\dot{\rho},H_z,\dot{z}$ & $t$     	& ${i=90\degr}$  \\
    \textbf{Hill} & $r,u,\Omega,\dot{r},H,H_z$ & $u$  & ${i=0\degr, 180\degr}$ 		\\
    \textbf{EccHill} & $H,H_z,\hat{f},\hat{g},\Omega,u$ & $u$  & ${i=0\degr, 180\degr}$  	\\
    \textbf{Ideal} 	& $\lambda_1, \lambda_2, \lambda_3, \lambda_4, H, C, S, \theta$ 	& $\theta$ 	& $-$      \\
    \hline
    \end{tabular}%
  \label{tab:CoordinateComparison}%
\end{table*}%

\section{Differential Algebra Techniques}\label{Sec3}

In this section we give a brief introduction to the DA framework and its application to the automatic computation of high-order expansions of the solution of ordinary differential equations (ODEs) and parametric implicit equations (PIEs). The interested reader can find further details in \citet{Berz1999v2} and \citet{M2013Nonlinear}.

\subsection{Differential algebra framework}\label{3A}

DA enables the efficient computation of the derivatives of functions within a computer environment. This is achieved by substituting the classical implementation of real algebra with the implementation of a new algebra of Taylor polynomials. Similarly to algorithms for floating point arithmetic, various algorithms were introduced in DA to treat common elementary functions, to perform composition of functions, to invert them and to solve nonlinear systems explicitly \citep{Berz1999v2}. In addition to these basic algebraic operations, operations for differentiation and integration were introduced in the algebra to complete the differential algebraic structure of DA. As a result, any deterministic function $f$ of $v$ variables that is $\mathcal{C}^{k+1}$ differentiable in the domain of interest $[-1,1]^v$ (these properties are assumed to hold for any function dealt with in this work) can be expanded into its Taylor polynomial up to an arbitrary order $k$ with limited computational effort. The DA used for the computations in this work was implemented by Dinamica \citep{ICATT2016} in the software DA Computational Engine (DACE), including all core DA functionality and a C++ interface.

\subsection{High-order expansion of the flow of ODEs}\label{3B}

An important application of DA is the automatic computation of the high-order Taylor expansion of the solution of ODEs with respect to the initial conditions and/or any parameter of the dynamics \citep{di2008application, M2013Nonlinear}. This can be achieved by replacing the classical floating point operations of the numerical integration scheme, including the evaluation of the right hand side of the ODE, with the corresponding DA-based operations. In this way, starting from the DA representation of the initial condition $\mathbf{X}_0$, the DA-based ODE integration supplies the Taylor expansion of the flow in $\mathbf{X}_0$ at all the integration steps, up to any final time $t_f$. Any explicit ODE integration scheme can be adapted to work in the DA framework in a straightforward way. 
The numerical integrator used for building the maps in this paper is a DA implementation of the 8th-order variable-stepsize Runge-Kutta integrator (RK8(7)) by \citet{PrinceDormand1981} with an 8th-order solution for propagation and 7th-order solution used for step size control and using an absolute tolerance of $10^{-12}$. Besides, the floating-point number version of this integrator is used for normal numerical propagation. Moreover, before propagating, all variables are scaled to the same order of magnitude using the length, velocity and time scaling factors: $R_e$, $\sqrt{\mu/R_e}$ and $\sqrt{R_e^3/\mu}$, respectively. Finally, if not mentioned otherwise, 5th-order Taylor expansions are used for the high-order mapping.

\subsection{High-order expansion of the solution of parametric implicit equations}\label{Param}
Satisfying constraints, such as boundary conditions, often requires finding the solution of an implicit equation
\begin{equation} \label{impeqclassic}
\boldsymbol c(\boldsymbol{x}) = 0,
\end{equation}
with $\boldsymbol{c} : \mathbb{R}^n \rightarrow \mathbb{R}^n $. This equation can be solved numerically using established numerical techniques, e.g. Newton's method.

Now suppose the vector function $\boldsymbol{c}$ depends explicitly on a vector of parameters $\boldsymbol{p}$, which yields the PIE
\begin{equation} \label{impeq}
\boldsymbol c( \boldsymbol{x}, \boldsymbol{p}) = 0.
\end{equation}
The solution of Eq.\,\eqref{impeq} is the function $\boldsymbol{x} = \boldsymbol f (\boldsymbol{p})$ that solves the equality for any value of $\boldsymbol{p}$. 

DA techniques can effectively solve the previous problem by representing $\boldsymbol{f}(\boldsymbol{p})$ in terms of its Taylor expansion with respect to the reference parameters $\boldsymbol{p}_0$. This result is achieved by first computing the Taylor expansion of $\boldsymbol{c}$ with respect to the reference values $\boldsymbol{x}_0$ and $\boldsymbol{p}_0$, and then applying partial inversion techniques as detailed in the work by \citet{di2008application}. The final result is 
\begin{equation}\label{xexp}
\boldsymbol{{x}} =  \mathcal{T}^k_{\boldsymbol f} (\boldsymbol{p}),
\end{equation}
which is the $k$-th order Taylor expansion of the solution of the PIE of Eq.\,\eqref{impeq}. For every value of $\boldsymbol{p}$, the approximate solution of $\boldsymbol c( \boldsymbol{x}, \boldsymbol{p}) = 0$ can be easily computed by evaluating the Taylor polynomial \eqref{xexp}. 
The accuracy of the approximation depends on both the order of the Taylor expansion and the displacement of $\boldsymbol{x}$ from its reference value $\boldsymbol{x}_0$.

The capability of expanding the solution of PIEs is of key importance for this work, because it used in the computation of high-order expansions of Poincar\'e maps.

\section{High-order mapping}
The High-Order Transfer Map (HOTM) method was developed by \citet{berz1987} and first applied to propagate perturbed Keplerian motion by \citet{Wittig2015High}. The method exploits the quasi-periodicity of perturbed orbits to efficiently map the orbital state over many revolutions and can be explained as follows. 

\begin{figure}
  \centering
	\includegraphics[width=0.5\columnwidth]{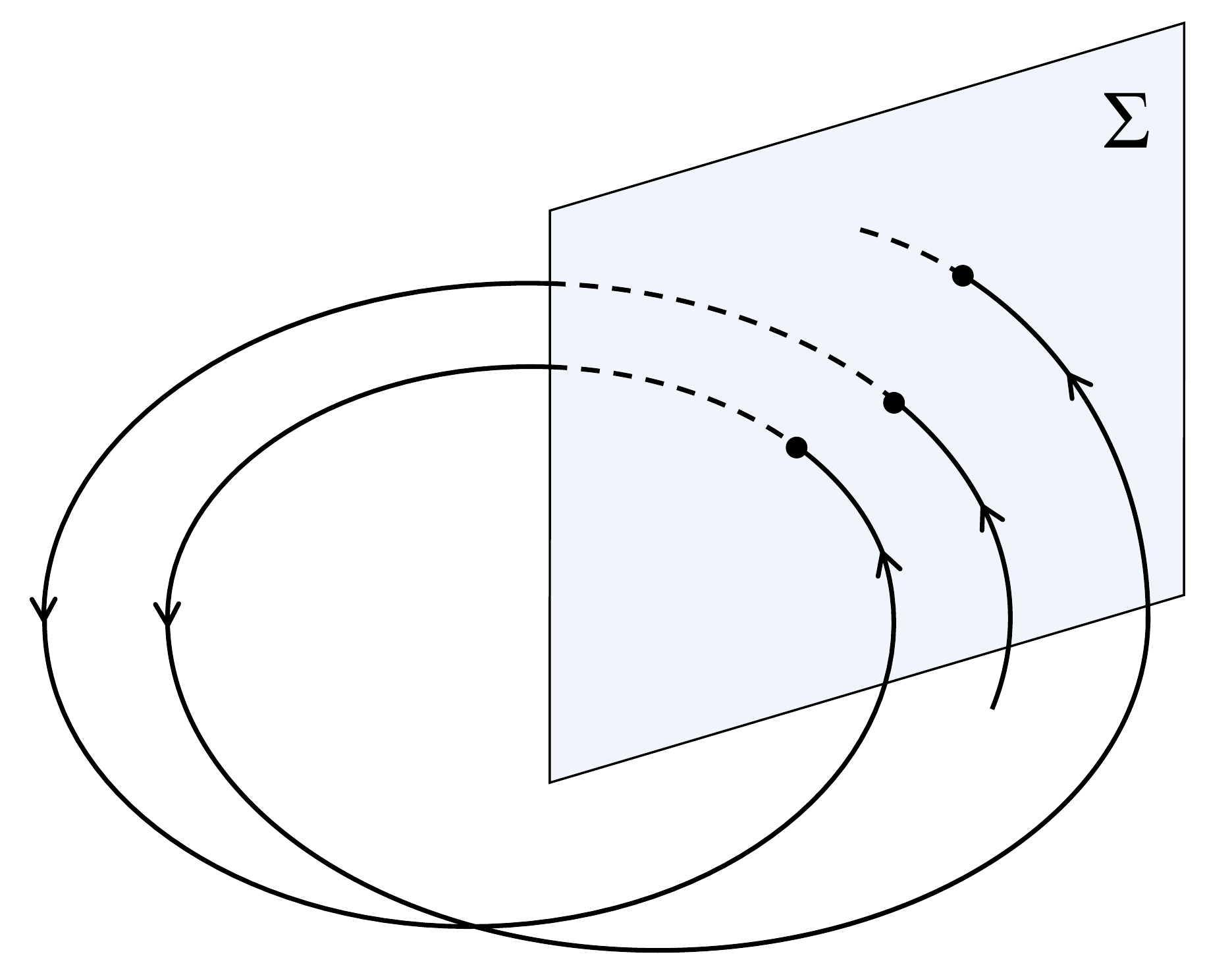}
    \caption{Orbits intersecting Poincar\'e section $\Sigma$.}
    \label{fig:PoincareSectionOrbits}
\end{figure}

Consider a hyperplane in the state space, that is the Poincar\'e section $\Sigma$, and quasi-periodic orbits that intersect the section, see Fig.~\ref{fig:PoincareSectionOrbits}. An orbit with its initial state $\mathbf{X}_0$ on the section, first leaves the section and then returns to it after one revolution. The function that maps the state of an orbit starting at $\Sigma$ over one revolution back onto $\Sigma$ is the Poincar\'e map $\Phi$. In other words, evaluating the map $\Phi$ for an orbital state $\mathbf{X}_0$ on the Poincar\'e section at the initial time $t_0$ gives the state of the orbit intersecting $\Sigma$ after one revolution, i.e. $\mathbf{X}_1 = \Phi(\mathbf{X}_0,t_0)$. This evaluation is usually time-consuming, because the dynamics are non-integrable and therefore has to be carried out numerically. 
However, since the dynamics are approximately periodic the state after one revolution $\mathbf{X}_1$ is close to the initial state $\mathbf{X}_0$ and in general two consecutive states $\mathbf{X}_{n+1}$ and $\mathbf{X}_n$ are close to each other. Therefore, by computing the Taylor expansion of the map around $\mathbf{X}_0$ and $t_0$ one can obtain a high-order approximation of $\Phi$ in a region close to $\mathbf{X}_0$ and at times close to $t_0$.

The HOTM method applies this idea by using DA to automatically compute a high-order Taylor expansion of $\Phi(\mathbf{X},t)$ around the expansion point $(\mathbf{X}_0,t_0)$ \citep{Wittig2015High}. This high-order map $\mathcal{T}_{\Phi}(\mathbf{X}, t)$ is an accurate approximation of $\Phi(\mathbf{X},t)$ for states that are close to $\mathbf{X}_0$ and times close to $t_0$. Besides, if the dynamics are autonomous, i.e. independent of time, then $\Phi$ only depends on the state, $\Phi=\Phi(\mathbf{X})$, and so does the high-order map, $\mathcal{T}_{\Phi}(\mathbf{X})$. In the following, we will assume that the dynamics are autonomous.

Mapping an orbit using the high-order map $\mathcal{T}_{\Phi}(\mathbf{X})$ is very fast because it only requires evaluating Taylor polynomials. In this way, the computationally expensive numerical integration of the continuous dynamics to evaluate $\Phi$ is replaced by efficient evaluation of the high-order map $\mathcal{T}_{\Phi}(\mathbf{X})$.

\subsection{High-order map computation}
\label{sec:mapcomputation}
To compute the high-order map $\mathcal{T}_{\Phi}(\mathbf{X})$ the orbit is propagated for one revolution in Taylor Differential Algebra with the state $\mathbf{X}_0$ initialized as a DA variable. Because perturbed Keplerian dynamics are only approximately periodic, the meaning of a revolution is ambiguous and any suitable definition of a revolution may be used. If the revolution is defined in space as the path between two consecutive crossings of a hyperplane, that is the Poincar\'e section $\Sigma$, then the map is a Poincar\'e map, see Fig.~\ref{fig:PoincareSectionOrbits}. Alternatively, a revolution can be defined by one period of the dynamics such that mapping occurs with the frequency of the dynamics and the map is a stroboscopic map.

\begin{figure}
  \centering
	\includegraphics[width=0.35\columnwidth]{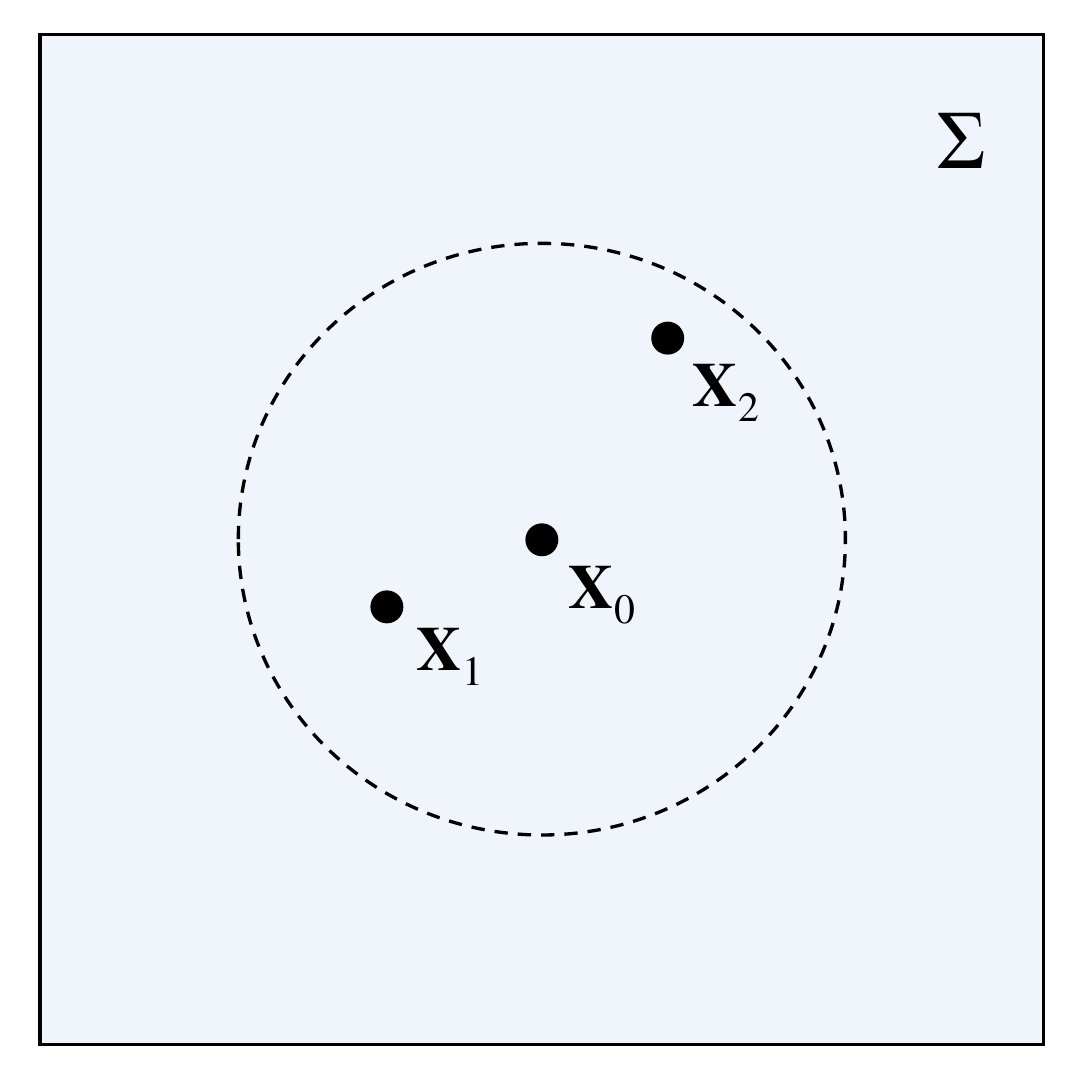}
    \caption{The Poincar\'e or stroboscopic map $\Phi$ is the map that maps the initial state $\mathbf{X}_0$ over one revolution to $\mathbf{X}_1$. The high-order map $\mathcal{T}_{\Phi}(\mathbf{X})$ is the high-order Taylor expansion of $\Phi$ around the expansion point $\mathbf{X}_0$. $\mathcal{T}_{\Phi}(\mathbf{X})$ can be used to accurately map $\mathbf{X}_1$ onto $\mathbf{X}_2$, and $\mathbf{X}_n$ onto $\mathbf{X}_{n+1}$, as long as the state $\mathbf{X}_n$ is in the domain close to the expansion point $\mathbf{X}_0$ where the truncation error is small. This accuracy domain is indicated by the dashed circle \citep{Wittig2015High}.}
    \label{fig:PoincareSectionHOTM}
\end{figure}

\subsubsection{Stroboscopic map computation}
For most element sets it is convenient to use stroboscopic mapping because the element sets contain a fast angular variable such that one revolution is defined as a change of $2\pi$ in the fast variable. In other words, the state is mapped at a fixed value of the fast angle\footnote{This way of stroboscopic mapping could be considered as Poincar\'e mapping using a Poincar\'e section that moves in inertial space but is fixed at a constant value of the fast angular coordinate.}. This approach is used for the COE, MEE, Hill, EccHill and Ideal element sets. To compute the map the fast angle is used as independent variable to DA-integrate the dynamics over $2\pi$ in the fast angle, see Table~\ref{tab:CoordinateComparison} for the used independent variables. For this purpose, the equations of motion are multiplied by the derivative of the time with respect to the fast variable $dt/ds$, where $s$ indicates the fast variable. One then obtains the equations of motion with respect to the fast variable and one equation for the evolution of time:

\begin{align}
	\frac{d\alpha_i}{ds} &= \frac{d\alpha_i}{dt} \frac{dt}{ds} ,  \\
	\frac{dt}{ds} &= \frac{1}{ds/dt} ,
\end{align} \\
where $\alpha_i$ are the orbital elements except the fast one.
The result of this DA propagation is straightforwardly the high-order stroboscopic map $\mathcal{T}_{\Phi}(\mathbf{X})$.

This approach works well for all orbital element sets except for the cylindrical coordinates for which the quasi-fast variable $\varphi$ is not defined at the poles.

\subsubsection{Poincar\'e map computation}
\label{sec:poincaremapcomp}
For the Cyl and CylHz element sets we use Poincar\'e mapping on the equatorial plane, i.e. at $z=0$. To achieve this, the nodal period (i.e., the time between two passages through the ascending node) is computed first and then $\mathcal{T}_{\Phi}$ is computed over one nodal period using time as independent variable. Because the nodal period is not constant and depends on the state, $T_n=T_n(\mathbf{X})$, it is approximated by a high-order Taylor expansion with respect to the initial state. The computation of the high-order Poincar\'e map for Cyl and CylHz coordinates is carried out as follows. 

First, the nodal period $T_{n}$ for the initial state $\mathbf{X}_0$, with $z_0 =0$, is computed numerically using the Keplerian orbital period $T = 2\pi\sqrt{a^3/\mu}$ as first guess. Then, both the nodal period and the initial state are initialized as DA variable and the dynamics are propagated over one nodal period in the DA framework, delivering 
\begin{equation}\label{prop}
\mathbf{X}_f = \mathcal{T}_{\mathbf{X}_f}(\mathbf{X},T_n).
\end{equation}
This state $\mathbf{X}_f$ must be on the equatorial plane, i.e. $z_f=0$, for any initial state $\mathbf{X}$ on $\Sigma$ to ensure Poincar\'e mapping. This condition is satisfied if $T_n$ is the nodal period corresponding to $\mathbf{X}$. Therefore, the computation of the Taylor approximation of the function $T_n(\mathbf{X})$ requires the solution of the PIE 
\begin{equation}\label{PIE1}
{z}_f = \mathcal{T}_{z_f}(\mathbf{X},T_n) = 0 ,
\end{equation}
in which $\mathcal{T}_{z_f}(\mathbf{X},T_n)$ is extracted from the map \eqref{prop}. This PIE, in which $T_n$ plays the role of the variable and $\mathbf{X}$ of the vector of parameters, is solved using the algorithm presented in Sec.~\ref{Param}, providing 
\begin{equation}\label{eq:TnExpansion}
{T}_n = \mathcal{T}_{T_n}(\mathbf{X}) ,
\end{equation}
i.e. the high-order Taylor approximation of the nodal period with respect to the initial state.

Finally, the high-order Taylor expansion of the Poincar\'e map is calculated by recomputing the Taylor expansion for $\mathbf{X}_f$ with ${T}_n$ initialized as $\mathcal{T}_{T_n}(\mathbf{X})$, so we get:
\begin{equation}\label{Poincare}
\mathcal{T}_{\Phi}(\mathbf{X}) = \mathcal{T}_{\mathbf{X}_f}(\mathbf{X},T_n(\mathbf{X})) = \mathcal{T}_{\mathbf{X}_f}(\mathbf{X}). 
\end{equation}
This map can be used to map a point on the Poincar\'e section $\Sigma$ at $z=0$ into its successive passage through the surface of section, and the expansion for $T_n$ \eqref{eq:TnExpansion} is used to keep track of the time of the passages.   

This approach of first solving for the nodal period and then building the high-order Poincar\'e map can be used for any element set. However, because it requires computing two high-order maps, which is time consuming, stroboscopic mapping is in general preferred.
Besides, it should be noted that using the argument of latitude $u$ as independent variable also enables Poincar\'e mapping on the equatorial plane, that is at $u=0$.

\subsection{Accuracy}
\label{sec:accuracy}
$\mathcal{T}_{\Phi}(\mathbf{X})$ is an approximation of the true transfer map $\Phi$ for states close to the expansion point $\mathbf{X}_0$. Therefore, the accuracy of the HOTM degrades when the state $\mathbf{X}_n$ drifts away from the initial state. If the dynamics also depend on time, then the accuracy also reduces as time passes. In this paper, we only consider autonomous perturbations, i.e. perturbations that do not explicitly depend on time.

In case of autonomous perturbations, the validity of the HOTM depends on the rate of change of the state and on the nonlinearity of the dynamics. If the dynamics are linear, then the HOTM is an exact approximation and is valid forever. If the dynamics are strongly nonlinear then the HOTM is only accurate for a small domain around the expansion point. The domain where the high-order map has a specific accuracy is called here the accuracy domain, see Fig.~\ref{fig:PoincareSectionHOTM}. This domain can be estimated using a method developed by \citet{wittig2015propagation} for estimating the truncation error of high-order Taylor expansions. By estimating the magnitude of higher-order terms, the distance from the expansion point where the Taylor series has a specific truncation error can be computed. For states inside the estimated domain the high-order map has a truncation error that is approximately smaller than the specified error.

\section{Fixed points of Poincar\'e maps}
\label{sec:fixedPoint}
The determination and the study of fixed points of Poincar\'e maps is one of the key topics in dynamical system theory. Moreover, fixed points are of great practical importance in astrodynamics because they provide ideal nominal orbits for space missions \citep{Coffey1994,dunham2003libration}. In particular, frozen orbits, i.e. orbits with stable eccentricity and argument of pericenter used by engineers since the early years of astrodynamics \citep{Coffey1994}, can be computed as fixed points of a reduced state in the zonal problem \citep{Broucke1994}. More recently, the centre manifold of these fixed points has been extensively studied by researchers with the aim of designing long-term and large amplitude relative bounded motion, suitable for formation flying missions \citep{koon2001j2,Xu2012,Baresi2017Design,baresi2017}. For these reasons the study of the motion of quasi-periodic orbits about a fixed point of the zonal problem is offered as a further test case in Sec.\,\ref{Tcases}.   

\subsection{Computation of fixed points}
Once the Taylor approximation of the Poincar\'e map is obtained as illustrated in Sec.~\ref{sec:mapcomputation}, the computation of its fixed points can be framed as a constraint satisfaction problem; i.e. find $\mathbf{X}^*$ such that 
\begin{equation}\label{eq:fixedPoint}
{\mathbf{X}^*} = \mathcal{T}_{\Phi}(\mathbf{X}^*). 
\end{equation}
This problem is here solved with the \textsc{matlab} nonlinear solver \texttt{fsolve}. Note that, as $\mathcal{T}_{\Phi}$ is a polynomial, the problem is reduced to finding the solution of a set of polynomial equations, for which all the derivatives required by the solver are readily available. In addition, in the zonal gravitational field the problem is reduced to a bidimensional constraint, because it is sufficient that distance $r$ and velocity $V_r=\dot{r}$ repeat after one nodal period. In EccHill elements this is equivalent to equal values of $\hat{f}$ and $\hat{g}$ when passing through the ascending node.

After a fixed point is computed, the behaviour of quasi-periodic orbits around it can be simply studied by the repetitive evaluation of $\mathcal{T}_{\Phi}$ for the set of initial conditions of interest, thus producing Poincar\'e section plots.

\section{Test Cases}\label{Tcases}
To test the high-order mapping using different coordinate sets, we focus on the main problem in artificial satellite theory, that is the $J_2$ perturbation only. The $J_2$ perturbation is the main perturbation in the low Earth orbit (LEO) region above 400 km, where other perturbations such as drag and third-body attraction are of second order \citep{wakker2015fundamentals}. This upper LEO region is the typical location for Earth observation satellites that often fly in ground-repeating and sun-synchronous orbits \citep{wakker2015fundamentals}. Therefore, we test the high-order mapping for orbits at 500 and 800 km altitude.

The effect of the Earth oblateness is characterized by secular changes in $\Omega$, $\omega$ and $M$ and short periodic changes in all orbital elements. The secular changes depend strongly on the inclination. Therefore, to investigate the characteristics of the use of different element sets, the mapping is carried out for orbits at different inclinations, namely $i$ is 0$\degr$, 30$\degr$, 63.4$\degr$ (i.e. the critical inclination) and 90$\degr$. 
At zero inclination the rate of change of $\Omega$ and $\omega$ is maximum. On the other hand, at the critical inclination the argument of perigee is frozen and whereas the ascending node does not precess when $i=90\degr$. Furthermore, to analyse the mapping performance for highly elliptical orbits (HEOs) we look at a Molniya-like orbit with an eccentricity of 0.74 at $i=30\degr$ and $i=63.4\degr$. 

Besides the $J_2$ perturbation, the orbital evolution of LEO orbits is mainly affected by higher-order zonal and drag perturbations. These perturbations affect the orbital period and the orientation of the orbital plane that are important for e.g. repeat ground track and sun-synchronous orbits. Therefore, the best-performing coordinate set is also tested for additional perturbations, namely $J_3$, $J_4$ and drag. 

The performances of the different element sets are tested by analysing the position error resulting from the high-order mapping. The position error can be computed by either including or neglecting the time of the mapped state. The position error including time is computed by comparing the mapped state with a numerically computed state that is propagated to the epoch of the mapped state in MEE, which is free of singularities. On the other hand, the error in position without considering time is calculated by comparing against numerical propagation using the same elements and same independent variable as the mapped state. For the Cyl and CylHz coordinates, the state on the Poincar\'e section is computed numerically to obtain the position error without time.

Finally, to demonstrate the potential of the method, we use high-order mapping to investigate the quasi-periodic orbits around a fixed point under $J_2$-$J_4$ perturbations by first computing the fixed point and subsequently mapping the orbits around it.

The initial conditions for all test cases and the initial guess for the fixed point are shown in Table~\ref{tab:testcases}. Notice that all orbits start on the equatorial plane.
The values of gravitational coefficients of the Earth and the drag parameters used for propagation are given in Table~\ref{tab:constants}.

\begin{table*}
  \centering
  \caption{Overview of test cases and initial osculating orbital elements. The initial state is always on the equatorial plane.}
    \begin{tabular}{llcccccc}
    \hline
    $\#$  & Test case & $a$ [km] & $e$ [-] & $i$ [$\degr$] & $\Omega$ [$\degr$] & $\omega$ [$\degr$] & $\nu$ [$\degr$] \\
    \hline
    1     & LEO, $J_2$ only & 6878.1363 & 0.01  & 30    & 30    & 30    & 330 \\
    2     & LEO, $J_2$ only & 6878.1363 & 0.01  & 0.0     & 30    & 30    & 330 \\
    3     & LEO, $J_2$ only & 6878.1363 & 0.01  & 63.4499 & 30    & 30    & 330 \\
    4     & LEO, $J_2$ only & 6878.1363 & 0.01  & 90    & 30    & 30    & 330 \\
    5     & HEO, $J_2$ only & 26561.7438 & 0.7411188 & 63.4428 & 30    & 270    & 90 \\
    6     & HEO, $J_2$ only & 26561.7438 & 0.7411188 & 30    & 30    & 270    & 90 \\
    7     & LEO, $J_2$-$J_4$ & 7178.1363 & 0.001 & 30    & 30    & 30    & 330 \\
    8     & LEO, $J_2$-$J_4$, drag & 6878.1363 & 0.01  & 30    & 30    & 30    & 330 \\
    9     & Fixed point, $J_2$-$J_4$ & 6878.1363 & 0.0 & 97.42 & 0.0     & 0.0    & 0.0 \\
    \hline
    \end{tabular}%
  \label{tab:testcases}%
\end{table*}%

\begin{table}
  \centering
  \caption{Values of constants and parameters used in the dynamical model.}
    \begin{tabular}{cc}
    \hline
    Parameter & Value \\
    \hline
    $\mu$ & 398600.4415 km$^3\,$s$^{-2}$ \\
    $R_e$ & 6378.1363 km \\
    $J_2$ & 0.001082626 \\
    $J_3$ & $-2.532411\times10^{-6}$ \\
    $J_4$ & $-1.619898\times10^{-6}$ \\
    $C_d$ & 2.2 \\
    $A/m$ & 0.0094736 m$^2\,$kg$^{-1}$ \\
    \hline
    \end{tabular}%
  \label{tab:constants}%
\end{table}%

\section{Results}
In this section, the results of high-order mapping using different element sets are presented and discussed. 
First, the results for a low-Earth low-eccentricity orbit and a highly-elliptical orbit under $J_2$ perturbation at different inclinations are discussed. After that, the mapping method is tested for higher zonal and drag perturbations and finally it is applied to investigate orbits around a fixed point.

\subsection{Low-Earth orbit}
In the following, the accuracy of high-order mapping of LEO orbits at different inclinations is analysed. For test case 1, first the performances of the established element sets are analysed and based on that the novel element set is discussed.

\subsubsection{Test case 1}
\begin{figure}
  \centering
	\includegraphics[width=0.8\columnwidth]{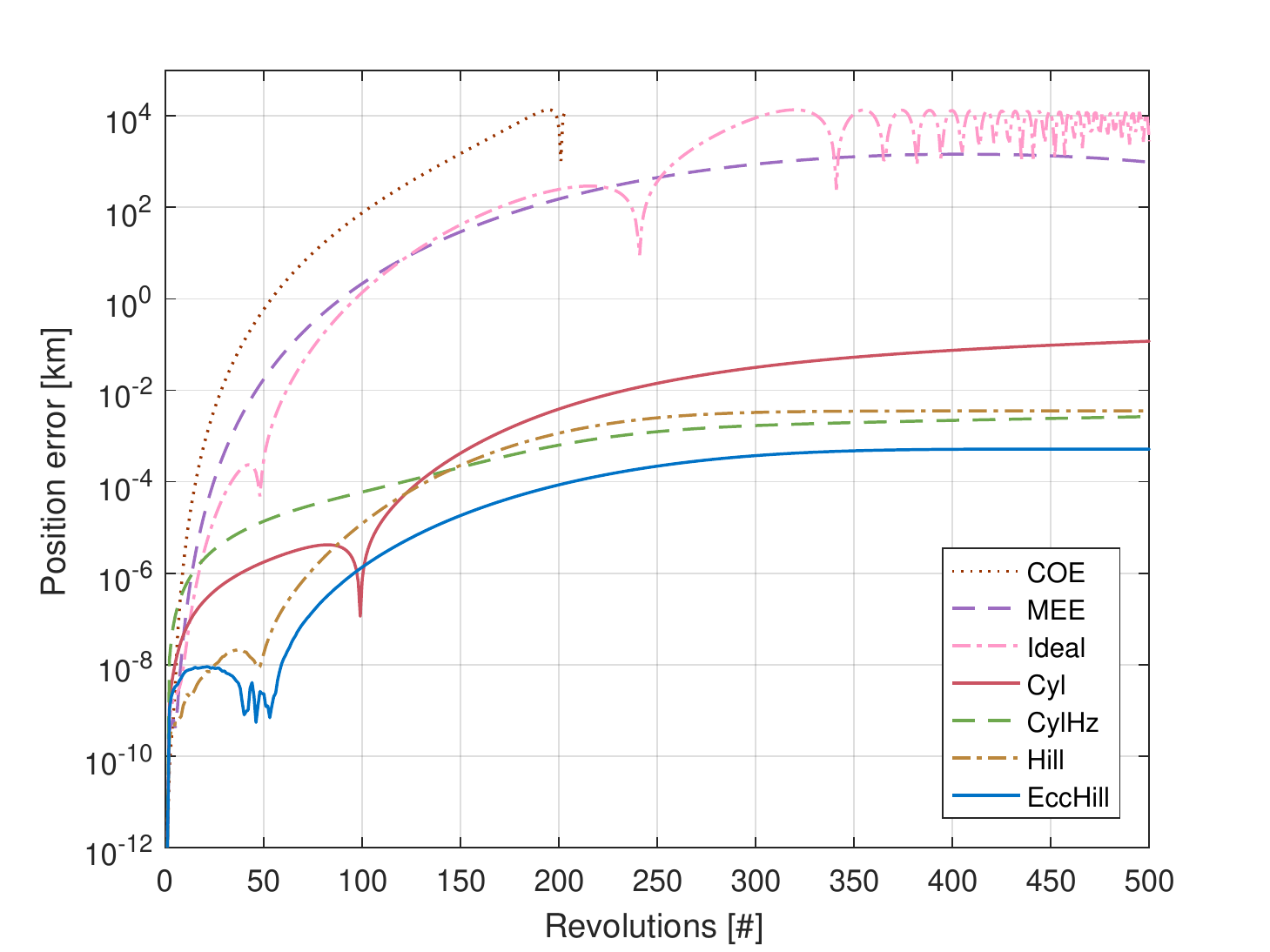}
    \caption{Position error for different element sets for test case 1 (LEO, $i=30\degr$, $J_2$ only) for 500 revolutions (33 days).}
    \label{fig:500km_e01_i30_500rev_PosError}
\end{figure}

\begin{figure}
  \centering
	\includegraphics[width=0.8\columnwidth]{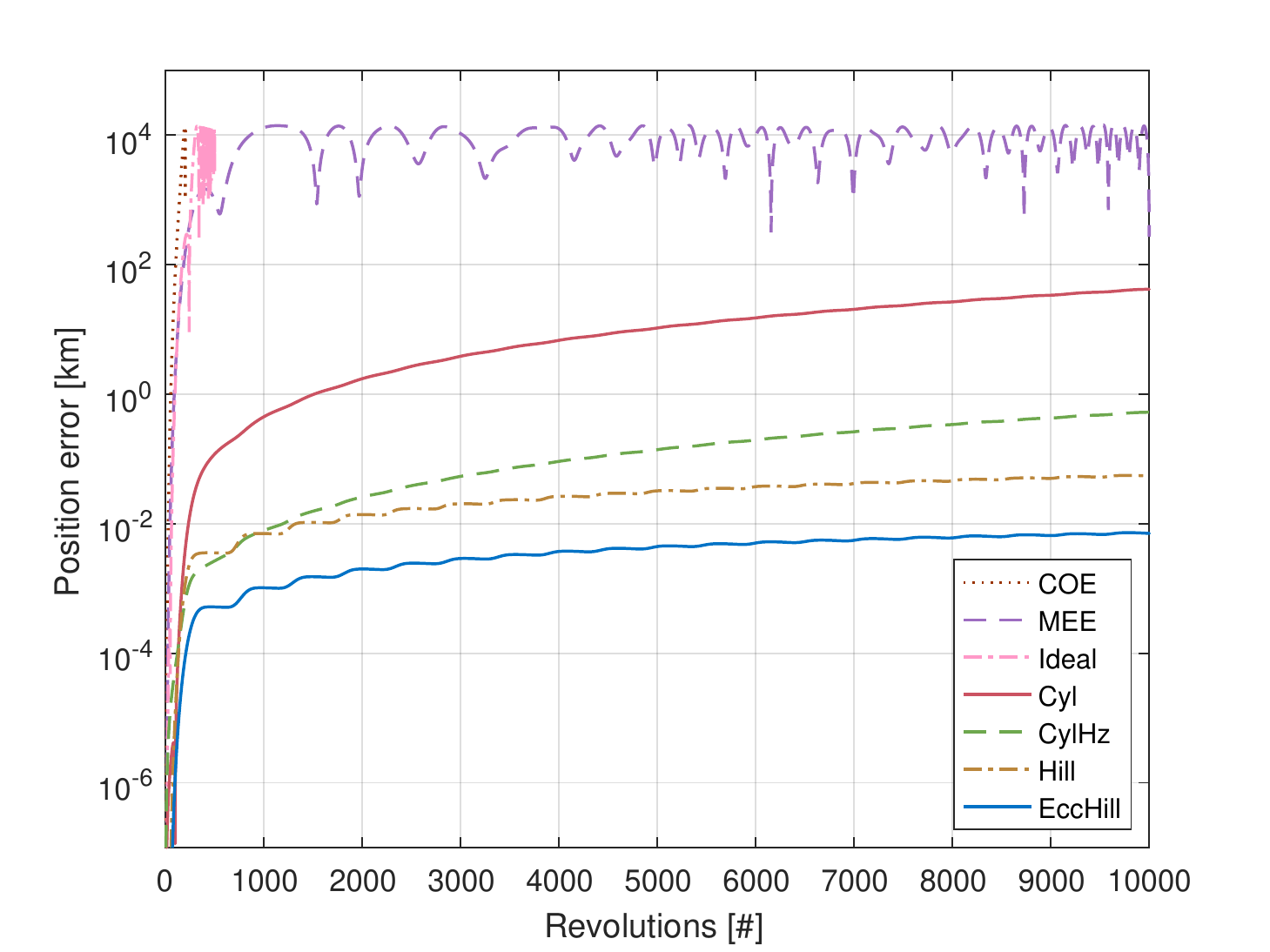}
    \caption{Position error for different element sets for test case 1 (LEO, $i=30\degr$, $J_2$ only) for 10,000 revolutions (655 days).}
    \label{fig:500km_e01_i30_10000rev_PosError}
\end{figure}

Figures \ref{fig:500km_e01_i30_500rev_PosError} and \ref{fig:500km_e01_i30_10000rev_PosError} show the position error of high-order mapping using different element sets for test case 1 ($i=30\degr$) for 500 and 10,000 revolutions, respectively. The results show that the COE, Ideal and MEE sets perform worst, having a position error larger than 1 km within 100 mappings. The CylHz and Hill sets, on the other hand, do not exceed an 1 km error within 10,000 mappings. This shows the large impact that the choice of elements has on the accuracy of the high-order mapping.

To determine the cause for the differences in performance, we analysed which specific variable in the different element sets caused the error to grow most. For this, the domain in which the high-order map has a truncation error less than $10^{-9}$ was computed, see Section~\ref{sec:accuracy}. If the variables remain inside this accuracy domain, then the Taylor series truncation error is very small (approximately less than $10^{-9}$). If, however, a variable drifts outsides the domain the accuracy of the high-order Taylor map decreases and the position error grows.
The element that is first to leave the accuracy domain is thus the main cause of large errors in the mapping.

Table~\ref{tab:ErrorElements} shows the number of mappings after which an element drifted outside the accuracy domain and which element it was for the different sets. These results show that the error indeed grows fastest for the COE, Ideal and MEE sets as the estimated truncation error exceeded $10^{-9}$ after 17 or less mappings.
This decrease in accuracy is caused by the secular drift in $\Omega$ and $\omega$ due to $J_2$. For the COE set, $\omega$ appears in the equations of motion in sine and cosine terms via the argument of latitude $u$. The Taylor series of sine and cosine functions are only accurate in a small domain, because low-order polynomials cannot accurately approximate these functions. As a result, the accuracy of the expansions reduces quickly when the state drifts away from the expansion point. 

\begin{table}
  \centering
  \caption{Number of mappings after which the specified element drifts out of the accuracy domain, which is the domain where the high-order map has an estimated truncation error less than $10^{-9}$, for test case 1 (LEO) and test case 6 (HEO).}
    \begin{tabular}{lcccc}
    \hline
    Element & \multicolumn{2}{c}{LEO} & \multicolumn{2}{c}{HEO} \\
\cline{2-5} set  & \multicolumn{1}{c}{\# of} &  & \multicolumn{1}{c}{\# of} &  \\
    	& \multicolumn{1}{c}{mappings} & Element & \multicolumn{1}{c}{mappings} & Element \\
    \hline
    COE   & 2     & $\omega$ 		& 14    & $\omega$ \\
    Ideal & 6     & $\lambda_2$ 	& 3     & $C$ \\
    MEE   & 17    & $k$   			& 6     & $f$ \\
    Cyl   & 18    & $\dot{\varphi}$ & 2     & $\rho$ \\
    CylHz & 40    & $\rho$ 			& 2     & $\rho$ \\
    Hill & 57    & $r$   			& 2     & $r$ \\
    EccHill & 126   & $\hat{f}$ 		& 3     & $\hat{f}$ \\
    \hline
    \end{tabular}%
  \label{tab:ErrorElements}%
\end{table}

For the Ideal and MEE sets, on the other hand, it is the secular change in $\Omega$ that causes the rapid decrease in accuracy. The $J_2$ perturbation does not depend on $\Omega$; however, the variation in $\Omega$ affects the components of the Ideal quaternion $\boldsymbol{\lambda}$ and the MEE parameters $h$ and $k$. In these elements, $\Omega$ appears together with the inclination. Because the $J_2$ perturbation depends strongly on the inclination (see equations \eqref{eq:J2radial}-\eqref{eq:J2normal}), the high-order maps are sensitive to changes in $h$ and $k$, and $\boldsymbol{\lambda}$. Consequently, as $\Omega$ changes secularly, the accuracy of the maps in MEE and Ideal elements reduces quickly. 

The high-order maps in Hill, Cyl and CylHz elements, on the other hand, are not affected by the change in $\Omega$, because $\Omega$ and $\varphi$, respectively, do not appear in the equations of motion. 
However, the drift in $\omega$ causes changes in $\dot{\varphi}$, $\rho$ and $r$ that reduces the accuracy of the mappings in Cyl, CylHz and Hill elements, respectively, see Table~\ref{tab:ErrorElements}. Besides, because the CylHz element set does not contain $\dot{\varphi}$ as variable but instead includes $H_z$, which is constant under zonal perturbations, the accuracy reduces less quickly compared to mapping in Cyl coordinates.

Based on these observations, a new set of elements was developed, namely the eccentric Hill variables. This set is a modification of the Hill variables using the elements $H$, $H_z$, $\Omega$ and $u$ and replacing $r$ and $\dot{r}$ by $\hat{f}$ and $\hat{g}$, see Section~\ref{sec:axialnodal}. 
The element $r$ is the main cause of error growth in the Hill variables mapping. The new elements $\hat{f}$ and $\hat{g}$, on the other hand, are similar to $f$ and $g$ in the MEE set that cause less error growth. The elements $H_z$, $\Omega$ and $u$ were kept on purpose, because zonal perturbations do not depend on $\Omega$ and do not change the value of $H_z$ such that the accuracy of the high-order map is not affected by $\Omega$ and $H_z$. In addition, the element $u$ enables both stroboscopic mapping in general and Poincar\'e mapping on the equatorial plane. The drawback of using the element $\Omega$ is that it causes singularities at $i=0\degr$ and $i=180\degr$. However, eliminating this singularity requires coupling of $i$ and $\Omega$, like in the MEE set, which strongly reduces the accuracy of the mapping.

Figure \ref{fig:500km_e01_i30_10000rev_PosError} shows that the accuracy of the high-order mapping in eccentric Hill variables is better than all other element sets and the position error is less than 10 m for 10,000 revolutions. This performance in accuracy is supported by Table~\ref{tab:ErrorElements} that shows that the EccHill elements remain longest inside the domain with low truncation error. 

\begin{figure}
  \centering
	\includegraphics[width=0.8\columnwidth]{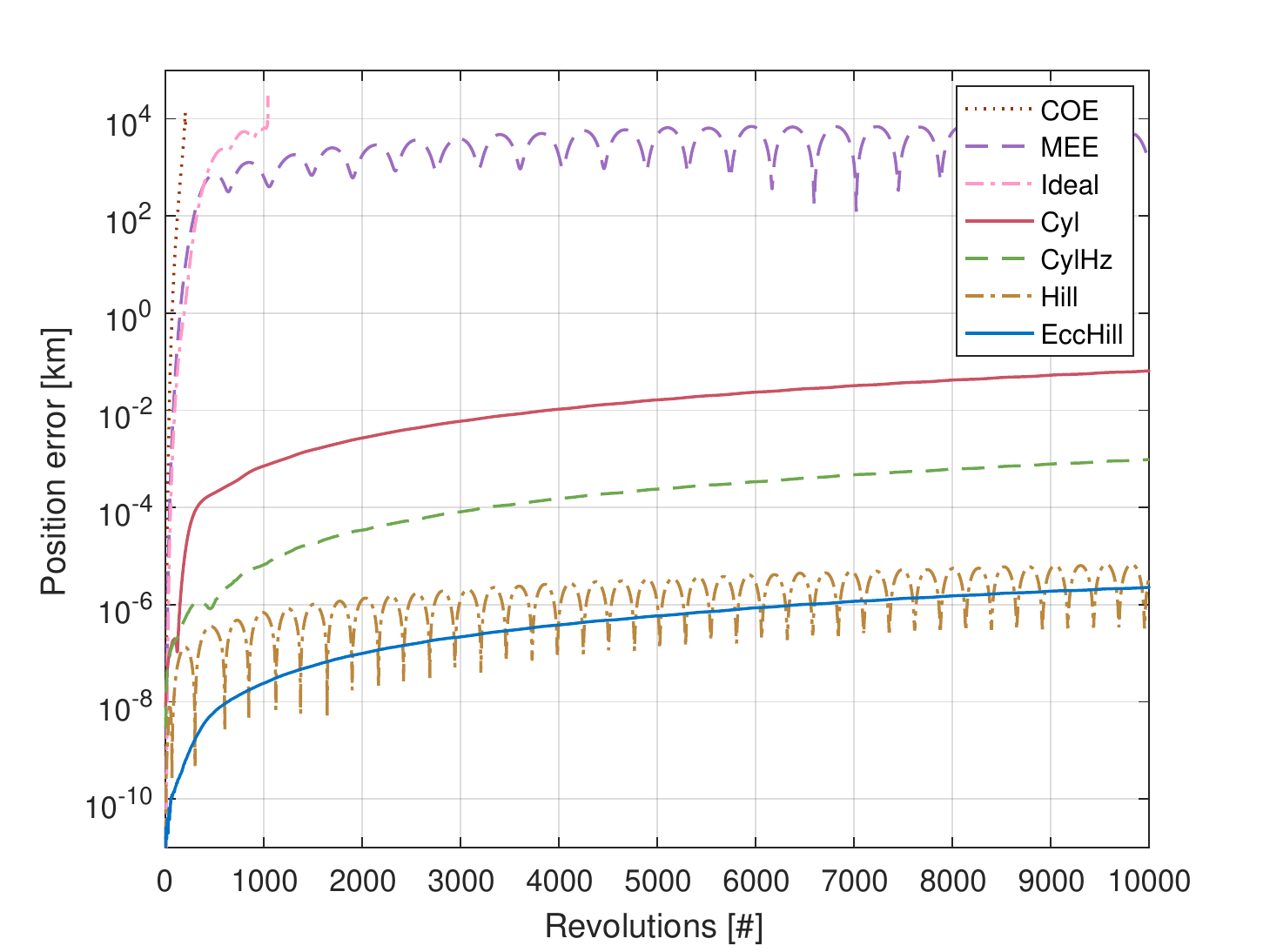}
    \caption{Position error without time error for different element sets for test case 1 (LEO, $i=30\degr$, $J_2$ only) for 10,000 revolutions (655 days).}
    \label{fig:500km_e01_i30_10000rev_PosError_ElementsOnly}
\end{figure}

The position errors shown in Figures \ref{fig:500km_e01_i30_500rev_PosError} and \ref{fig:500km_e01_i30_10000rev_PosError} can be decomposed into an error in the mapped elements and an error in the mapped time. For studying the evolution of an orbit the exact time of the mapping is not always of interest and only the elements need to be computed accurately. Figure~\ref{fig:500km_e01_i30_10000rev_PosError_ElementsOnly} shows the position error considering only the elements and not the time of the mappings. The errors are much smaller compared with Fig.~\ref{fig:500km_e01_i30_10000rev_PosError} which indicates that the time-included position error is for most part due to an error in time. Besides, the mapping errors using Hill and EccHill elements are extremely small with the maximum error lower than 1 cm (for 10,000 revolutions).

Finally, a key feature of the high-order map is that computationally expensive numerical propagation can be replaced by efficient evaluation of the map.
Table~\ref{tab:ComputationTimesLEOJ2i30} shows the CPU times required for mapping the LEO orbit of test case 1 using a high-order map and using numerical propagation. In addition, the CPU times for building the high-order map and for the actual mapping by evaluating the Taylor expansions are shown. First of all, it is clear that high-order mapping (including building the map) is an order of magnitude faster than numerical mapping for all elements sets. Secondly, computing the high-order map requires much more time than mapping the orbit by evaluating the Taylor map 10,000 times. This means that once a high-order map has been calculated, an orbit can be propagated very quickly for thousands of revolutions.

Considering the different element sets, mapping in Hill and EccHill elements is fastest, because their equations of motion for the $J_2$ perturbed problem are simple and thus fast to evaluate, and the propagation required few integration steps. In contrast, using Cyl, CylHz or Ideal elements for mapping demands the most computational effort. Cyl and CylHz coordinates require the calculation of two high-order maps to carry out Poincar\'e mapping (see Section~
\ref{sec:cylcoord}) and in the Ideal elements' dynamics many transformations between the Ideal and inertial frame are needed, which affects the computation time.

\begin{table}
  \centering
  \caption{CPU times in milliseconds for high-order mapping (building the map, mapping and total using 5th-order Taylor expansions) and numerical propagation of a LEO orbit (test case 1) for 10,000 revolutions. $^\ast$Includes computing additional high-order map to solve for Poincar\'e map condition, see Section~\ref{sec:poincaremapcomp}.}
    \begin{tabular}{lcccc}
    \hline
    Element set & \multicolumn{3}{c}{High-order mapping [ms]} & Numerical [ms] \\
          & Build map & Mapping & Total &  \\
    \hline
    COE   & 765   & 16    & 781   & 18299 \\
    MEE   & 983   & 16    & 999   & 9812 \\
    Hill & 140   & 16    & 156   & 5741 \\
    EccHill & 266   & 15    & 281   & 6162 \\
    Cyl   & 1376$^\ast$  & 31    & 1407  & 10561 \\
    CylHz & 1154$^\ast$  & 31    & 1185  & 9812 \\
    Ideal & 1919  & 93    & 2012  & 14321 \\
    \hline
    \end{tabular}%
  \label{tab:ComputationTimesLEOJ2i30}%
\end{table}%

\subsubsection{Test case 2}
At zero inclination and only $J_2$ perturbation, the transverse and normal perturbing forces, $f_t$ and $f_n$, are zero (see \eqref{eq:J2transverse} and \eqref{eq:J2normal}) and therefore $i$ and $H$ are constant. Furthermore, $\Omega$ and $\omega$ are not well defined at zero inclination and consequently the COE, Hill and EccHill sets become singular. However, in the special case of only even zonal harmonics, $f_n$ is zero which negates the singularity (e.g. $\sin{i}$ in the denominator of equations \eqref{eq:dRAANdt} and \eqref{eq:dargPerdt} is cancelled out by $\sin{i}$ in equation \eqref{eq:J2normal} for $f_n$) and the COE, Hill and EccHill elements can be used. For Cyl and CylHz coordinates, using the equatorial plane as Poincar\'e section is not possible when $i=0\degr$, therefore instead the Poincar\'e section is set at $\varphi=\varphi_0=30\degr$.

\begin{figure}
  \centering
	\includegraphics[width=0.8\columnwidth]{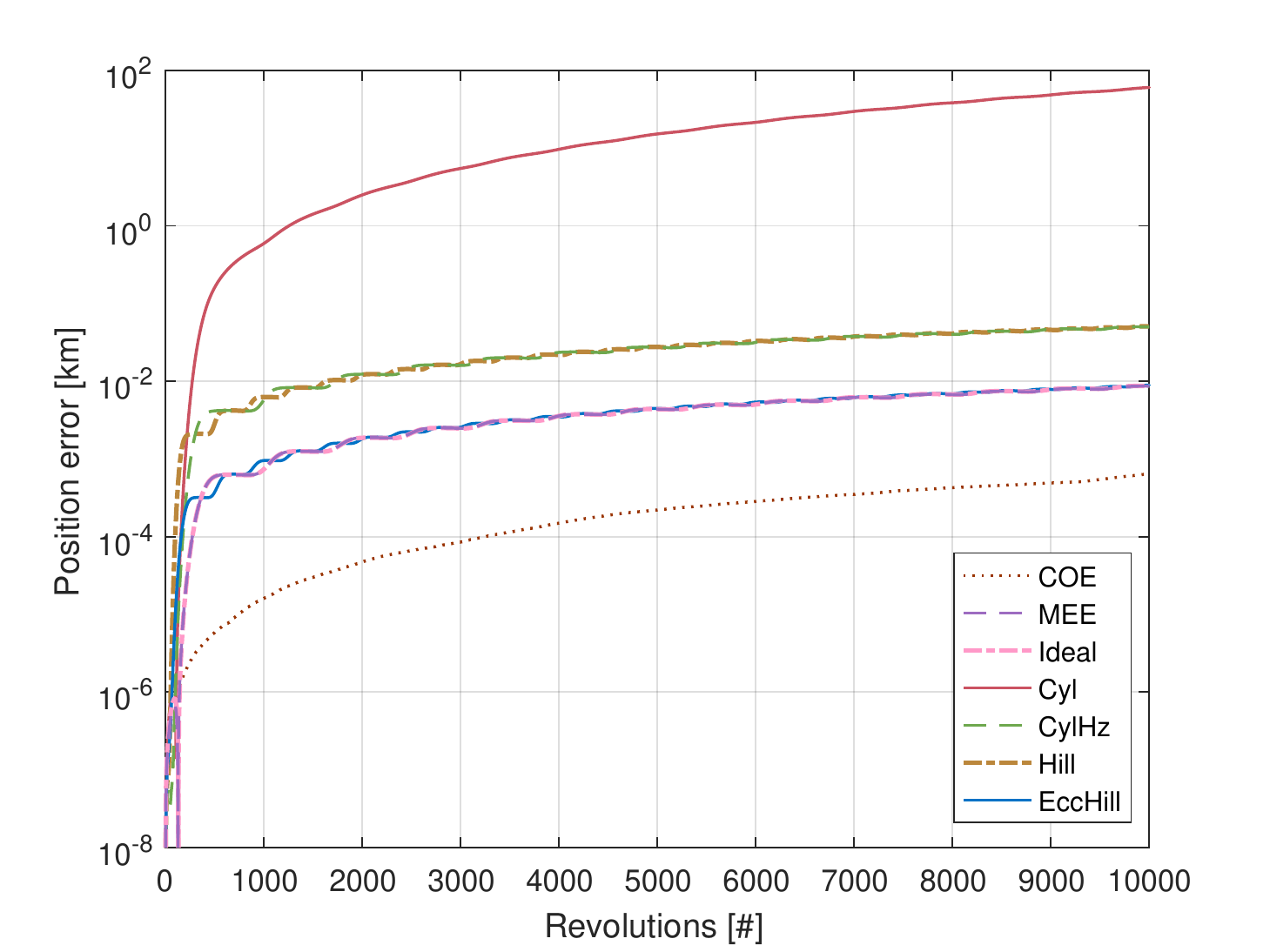}
    \caption{Position error for different element sets for test case 2 (LEO, $i=0\degr$, $J_2$ only) for 10,000 revolutions (655 days).}
    \label{fig:500km_e01_i0_10000rev_PosError}
\end{figure}

Figure~\ref{fig:500km_e01_i0_10000rev_PosError} shows the position error resulting from mapping a LEO orbit at $i=0\degr$ using different element sets. The COE set performs very well, because the drift in $\omega$ does not affect the mapping accuracy since $\omega$ is a longitudinal angle at $i=0\degr$ and $J_2$ does not depend on the longitude. The MEE, Ideal and EccHill elements perform equally well, because the MEE elements $h$ and $k$, the Ideal quaternion $\boldsymbol{\lambda}$ and EccHill variables $H$ and $H_z$ are constant and the other accuracy-affecting elements $f$ and $g$ (MEE), $C$ and $S$ (Ideal) and $\hat{f}$ and $\hat{g}$ (EccHill) evolve essentially the same at zero inclination. The Hill variables and CylHz set also perform alike because $r$ and $\dot{r}$, and $\rho$ and $\dot{\rho}$, respectively, vary similarly for $i=0\degr$. Finally, the Cyl coordinates achieve the lowest accuracy, but similar to test case 1 due to the change in $\dot{\varphi}$. 

\subsubsection{Test case 3}
\begin{figure}
  \centering
	\includegraphics[width=0.8\columnwidth]{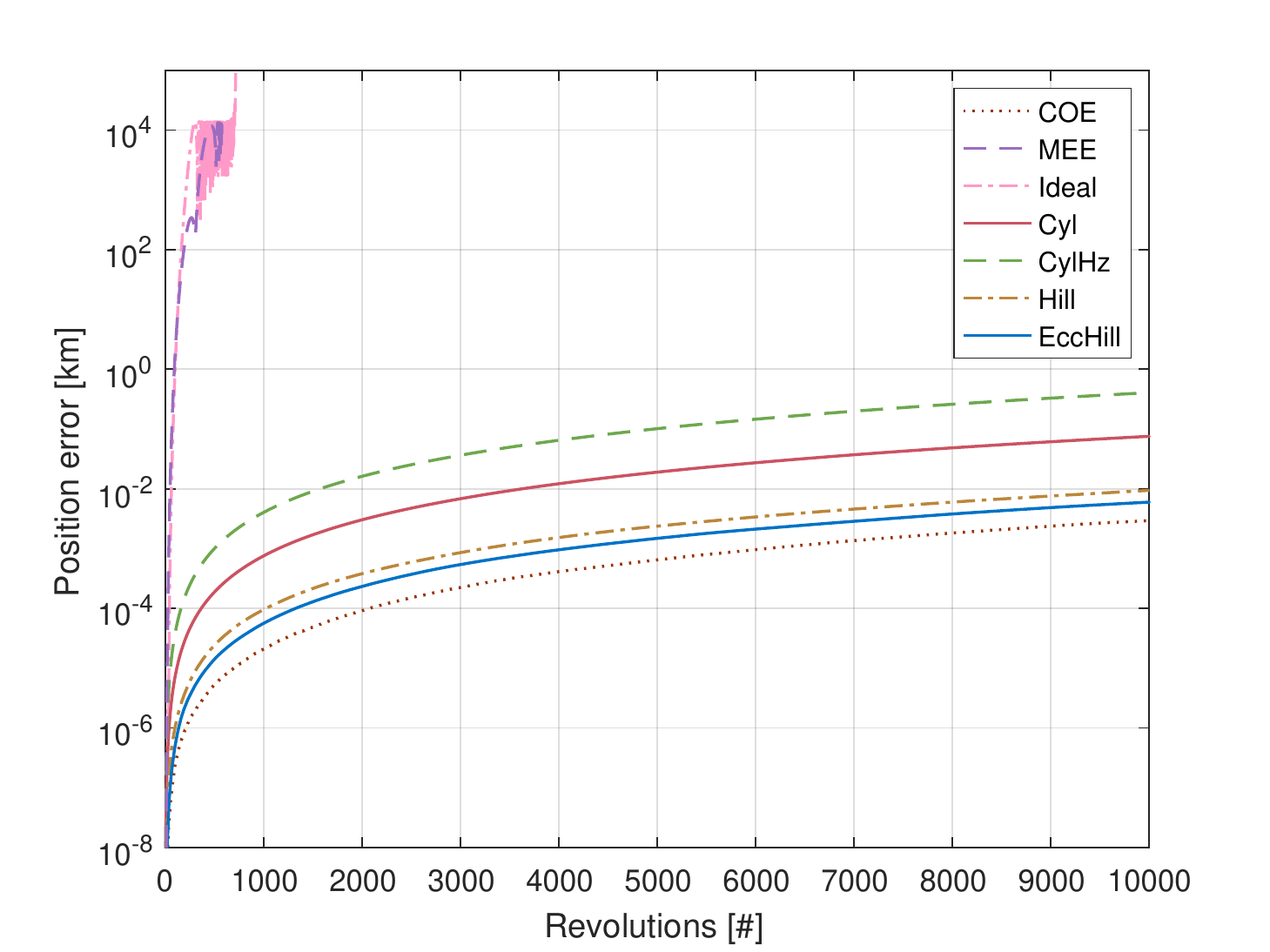}
    \caption{Position error for different element sets for test case 3 (LEO, $i=63.4\degr$, $J_2$ only) for 10,000 revolutions (655 days).}
    \label{fig:500km_e01_i63_10000rev_PosError}
\end{figure}

At an inclination of $63.4\degr$, called the critical inclination, the secular variation of $\omega$ is zero, i.e. $\omega$ is fixed. The accuracy of the high-order mapping for an LEO orbit at critical inclination is shown in Fig.~\ref{fig:500km_e01_i63_10000rev_PosError}. The MEE, Ideal, CylHz and EccHill elements perform similar to test case 1, because they do not profit from the zero variation in $\omega$. The COE, Cyl and Hill sets, on the other hand, achieve better accuracies. For COE and Hill variables this is expected because precession of the pericenter affects the mapping accuracy. The mapping using Cyl coordinates improves because $\dot{\varphi}$ does not vary since $\omega$ is fixed.

\subsubsection{Test case 4}
\begin{figure}
  \centering
	\includegraphics[width=0.8\columnwidth]{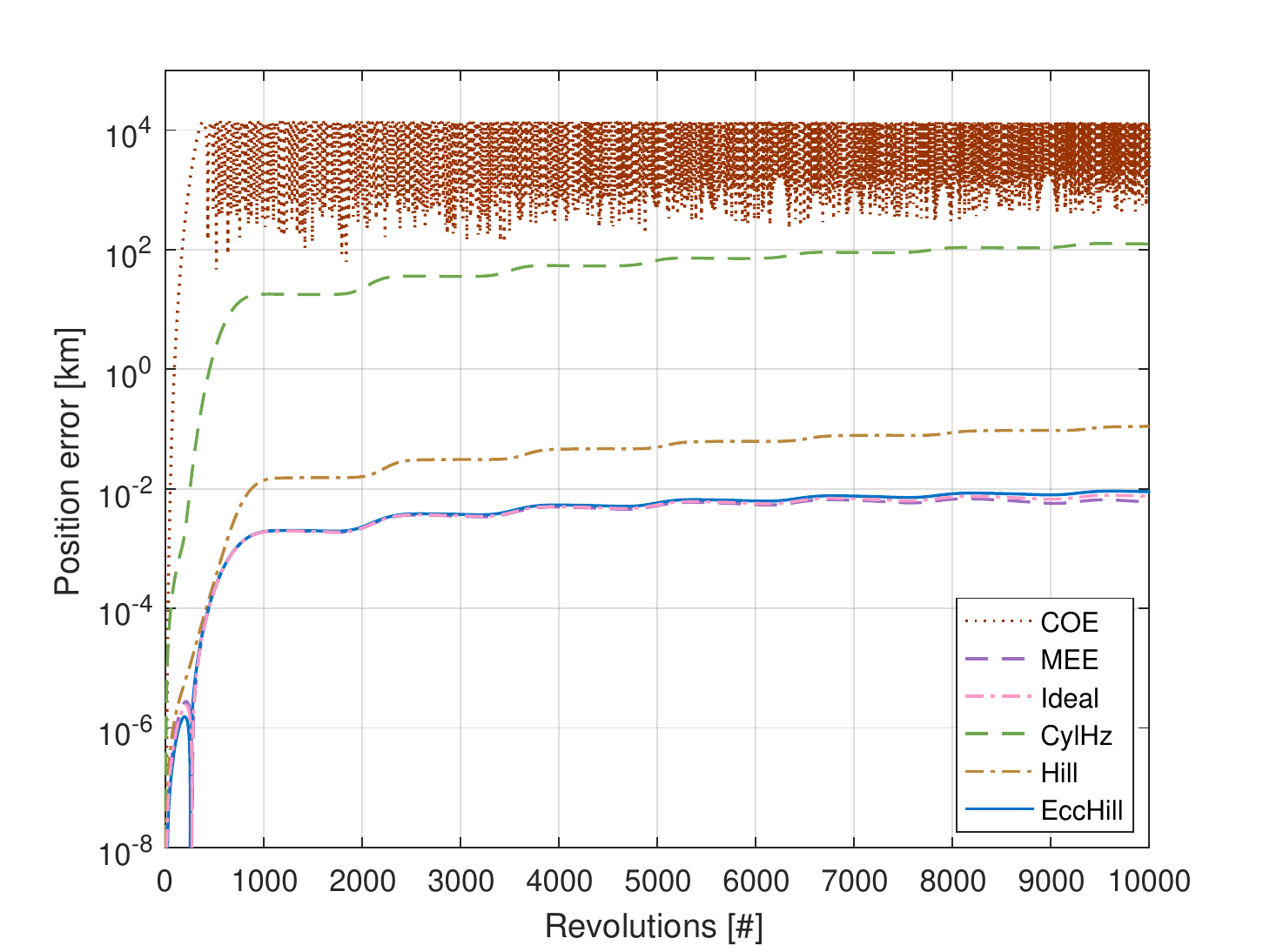}
    \caption{Position error for different element sets for test case 4 (LEO, $i=90\degr$, $J_2$ only) for 10,000 revolutions (655 days).}
    \label{fig:500km_e01_i90_10000rev_PosError}
\end{figure}

For polar orbits, i.e. $i=90\degr$, the ascending node $\Omega$ is fixed. The position errors for a LEO orbit at $i=90\degr$ are shown in Fig.~\ref{fig:500km_e01_i90_10000rev_PosError}. As in the zero inclination case, the MEE, Ideal and EccHill elements sets perform equally well, because the MEE elements $h$ and $k$ and the Ideal quaternion $\boldsymbol{\lambda}$ are constant (since $f_n=0$) and therefore do not affect the accuracy. The COE and CylHz set perform poorly due to the precession of the pericenter. The Cyl coordinates could not be used, because they are singular at the poles.
\\

\noindent Overall, the LEO orbit test results have shown that the newly introducted element set, the eccentric Hill variables, perform best. Only the COE set performed better for the special cases of zero and critical inclination. In addition, the main causes of error growth were found be the drift in $\omega$, $\Omega$ and $r$.

\subsection{Highly elliptical orbit}
The dynamics of HEO orbits is more nonlinear than the motion of LEO satellites because the variation of the perturbing forces and coordinates, such as $r$, over one orbital revolution is much larger. 
Therefore, we test the performance for an orbit with a high eccentricity, namely a Molniya type of orbit.
The initial orbital elements for the HEO test cases are given in Table~\ref{tab:testcases}. 

\subsubsection{Test case 5}
\begin{figure}
  \centering
	\includegraphics[width=0.8\columnwidth]{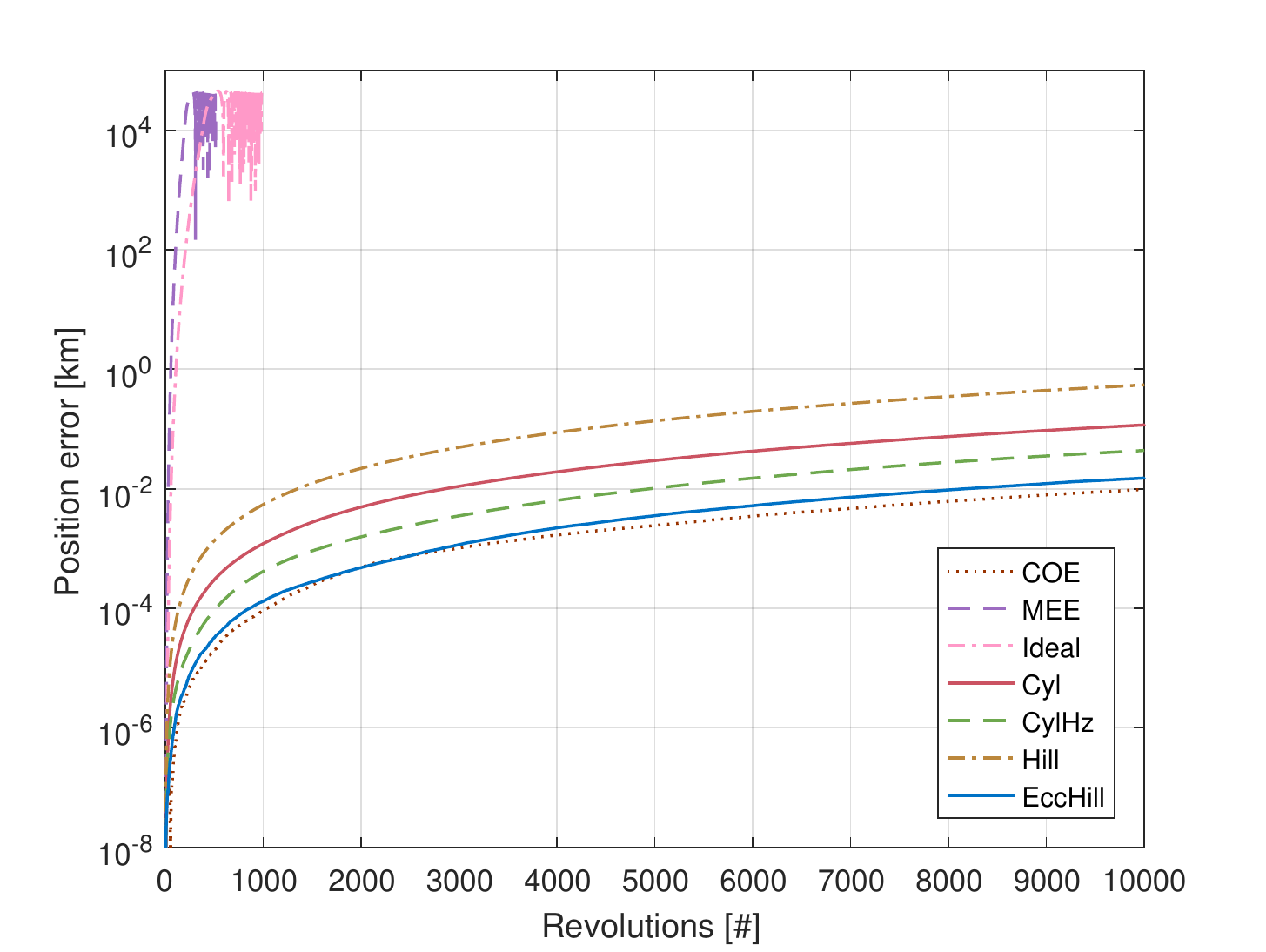}
    \caption{Position error for different element sets for test case 5 (HEO, $i=63.4\degr$, $J_2$ only) for 10,000 revolutions (13.6 years).}
    \label{fig:Molniya_i63_10000rev_PosError}
\end{figure}

If the inclination is set to the critical inclination such that the argument of pericenter is fixed, then the performance of the different element sets is similar to the LEO orbit case at critical inclination, see Fig.~\ref{fig:Molniya_i63_10000rev_PosError}. This means that the eccentricity does not affect the accuracy much when $\omega$ is constant. 

\subsubsection{Test case 6}
\begin{figure}
  \centering
	\includegraphics[width=0.8\columnwidth]{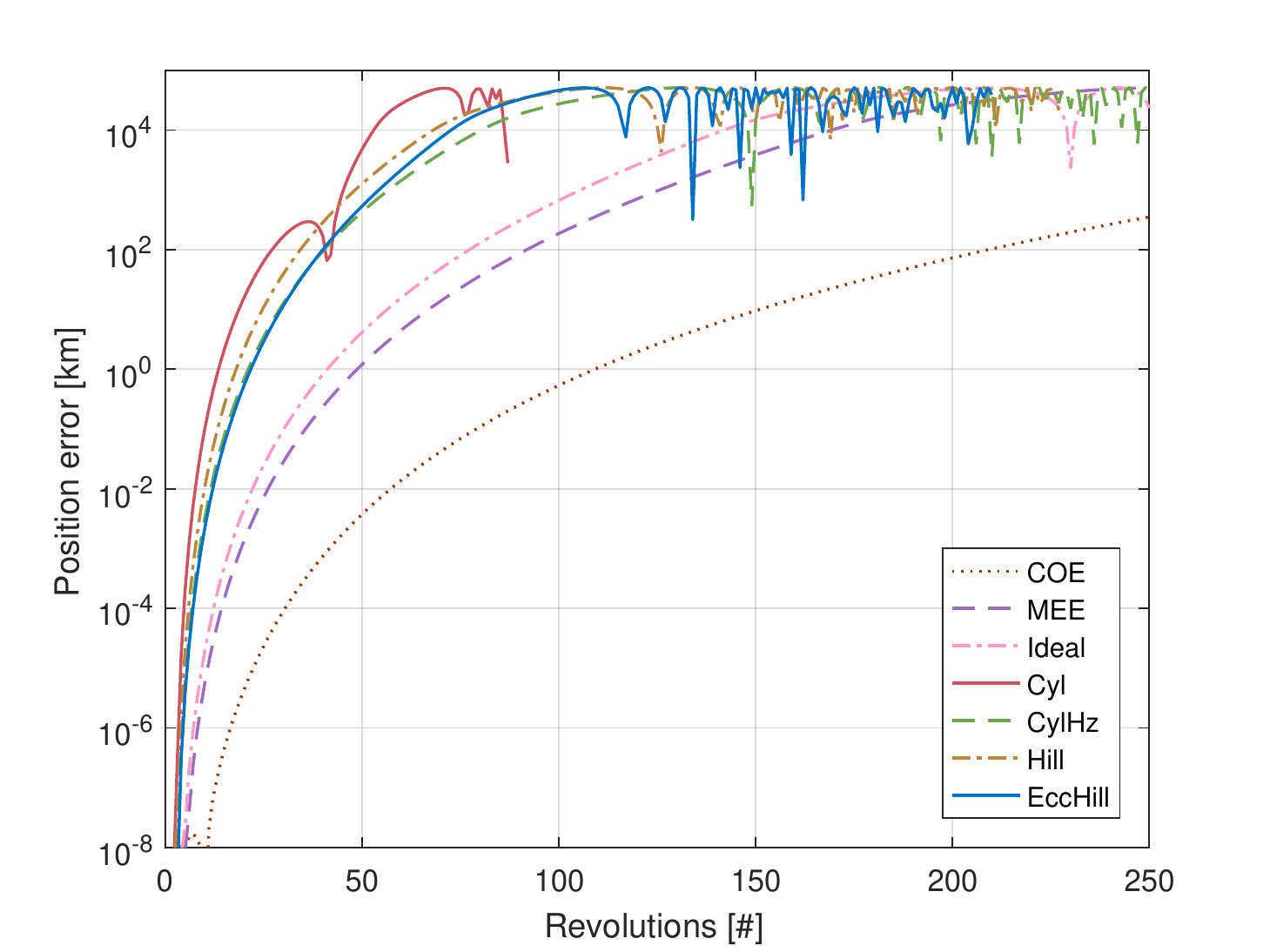}
    \caption{Position error for different element sets for test case 6 (HEO, $i=30\degr$, $J_2$ only) for 250 revolutions (125 days).}
    \label{fig:Molniya_i30_250rev_PosError}
\end{figure}

When the inclination is set to a different value, in this case 30$\degr$, then the performance of the high-order mapping compared with a near-circular orbit reduces strongly, see Fig.~\ref{fig:Molniya_i30_250rev_PosError}. The position error using Cyl, CylHz, Hill and EccHill elements grows larger than 1 km within 22 mappings. The COE, MEE and Ideal perform best, which is completely opposite to the LEO case, see Fig.~\ref{fig:500km_e01_i30_500rev_PosError}. The mapping in COE is most accurate and performs better than in the LEO case because $\omega$ precesses slower due to the larger semi-major axis. These results are supported by Table~\ref{tab:ErrorElements} that shows the amount of mappings after which the state drifts outside the accuracy domain. The table indicates the large decrease in accuracy as the number of mappings after which the truncation becomes larger than $10^{-9}$ is much lower than in the LEO for all element sets except for the COE set.
The Cyl, CylHz and Hill elements perform poorly because the mapped radial distance changes quickly as the elliptical orbit precesses. The variation in MEE variables $f$ and $g$, EccHill variables $\hat{f}$ and $\hat{g}$ and Ideal elements $C$ and $S$ is also larger compared the LEO case due to the increased eccentricity. However, the MEE variables $f$ and $g$ vary more slowly than $\hat{f}$, $\hat{g}$, $C$ and $S$, because the changes in $\Omega$ and $\omega$ are in opposite direction.

\begin{figure}
  \centering
	\includegraphics[width=0.8\columnwidth]{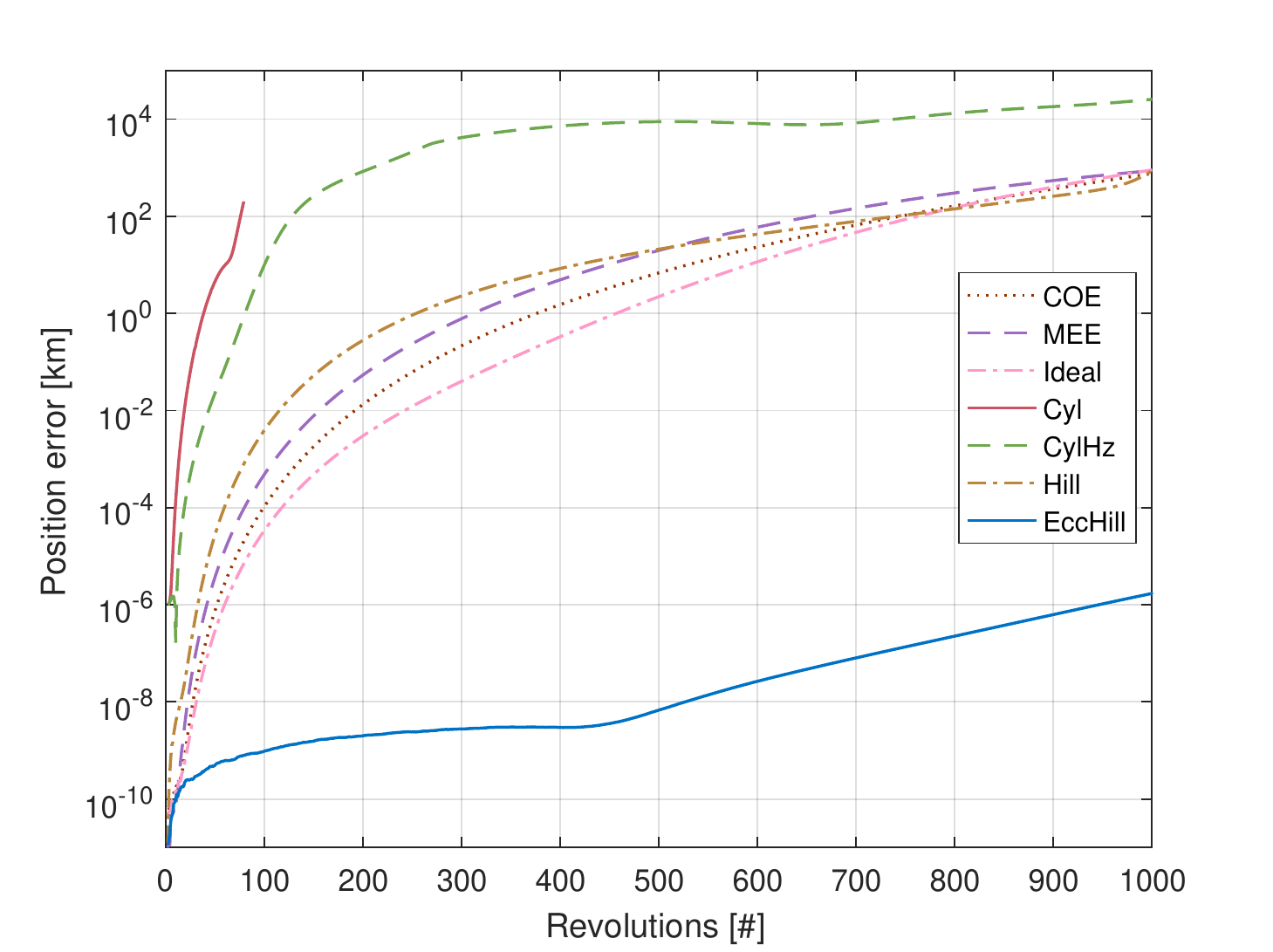}
    \caption{Position error without time error for different element sets for test case 6 (HEO, $i=30\degr$, $J_2$ only) for 1000 revolutions (498 days).}
    \label{fig:Molniya_i30_1000rev_PosError_ElementsOnly}
\end{figure}

Figure~\ref{fig:Molniya_i30_1000rev_PosError_ElementsOnly} shows position error when only considering the elements and neglecting time for the HEO orbit for 1000 revolutions. The results are much better and the number of revolutions before the position error grows larger than 1 km is more than 250 revolutions (125 days) for all element sets except the cylindrical coordinates. This means that the low accuracy compared to the LEO case is mainly caused by errors in the expansion of time and to a lesser extend by errors in the element expansions. For the Cyl and CylHz sets, time is the independent variable for the mapping and the mapped elements thus depend on time. Therefore, the position error is not much better when only the elements are considered and not the timing.

\begin{figure}
  \centering
	\includegraphics[width=0.8\columnwidth]{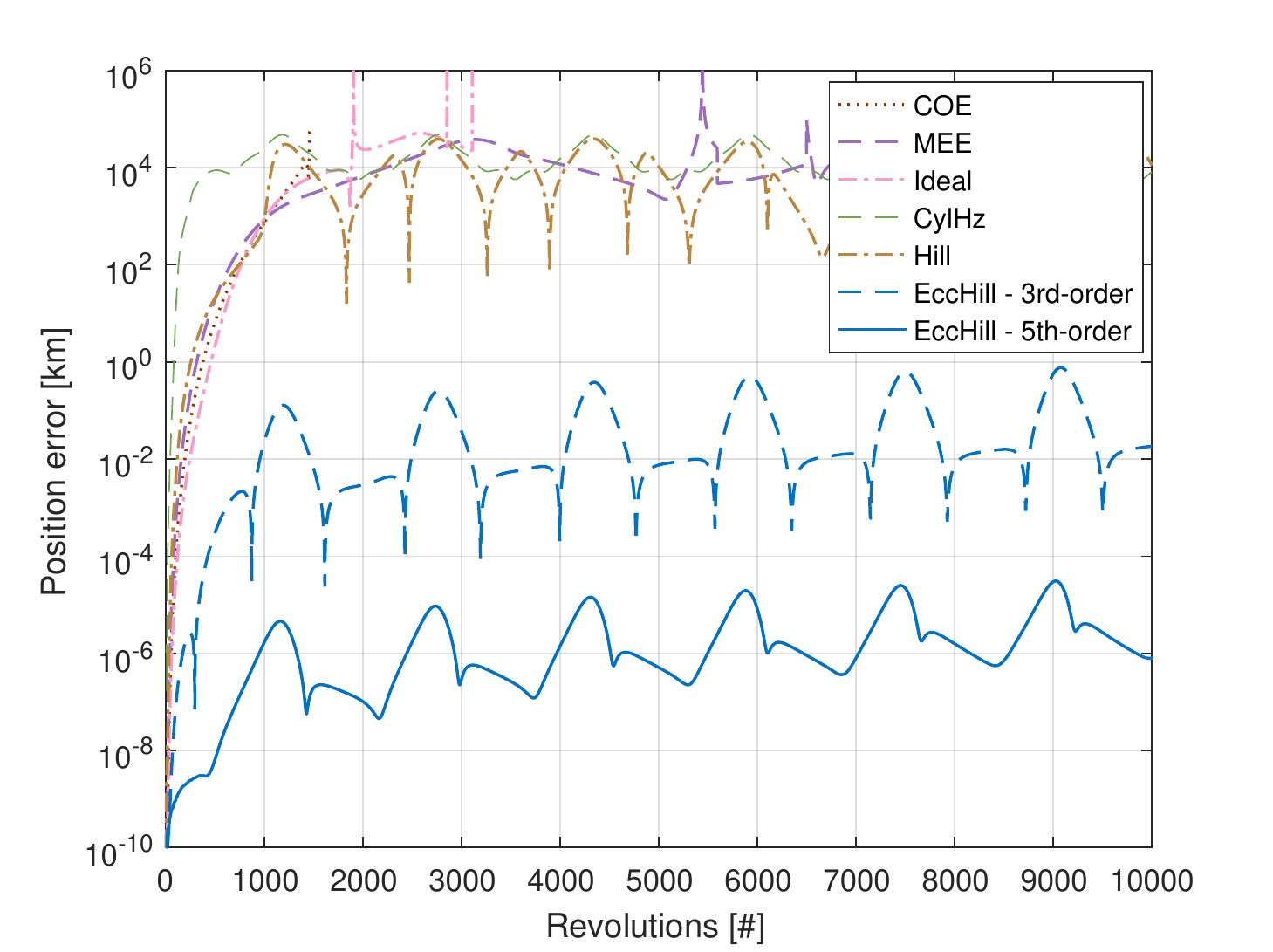}
    \caption{Position error without time error for different element sets for test case 6 (HEO, $i=30\degr$, $J_2$ only) for 10,000 revolutions (13.6 years).}
    \label{fig:Molniya_i30_10000rev_PosError_ElementsOnly}
\end{figure}

Figure~\ref{fig:Molniya_i30_10000rev_PosError_ElementsOnly} shows the element-only error for 10,000 revolutions. The EccHill element set performs extremely well with a position error less than 3 cm. The EccHill elements are almost mapped exactly and all error can thus be said to be caused by the expansion for time. The mapping of EccHill elements is so accurate that even using 3th-order Taylor expansions the element-only position error is less than 1 km for 10,000 revolutions. 
These accuracies are very high considering that mapping the Molniya orbit for 10,000 revolutions means propagating the osculating elements for more than 13 years (4980 days). 

The computation times for high-order mapping and numerical propagation of the HEO orbit are shown in Table~\ref{tab:ComputationTimesHEOJ2i30}. As in the LEO case, the high-order mapping including building the map is an order of magnitude faster than mapping using numerical propagation, and the mapping in Hill and EccHill elements is fastest. In general, the CPU times are about two times higher compared to the LEO orbit because of the nonlinearity of the HEO dynamics. Evaluating the high-order maps, on the other hand, requires the same amount of time and is about three orders of magnitude faster than numerical mapping.\\

\noindent The results for the LEO and HEO orbits have shown that the EccHill elements perform best in terms of accuracy. In the next section, the high-order mapping using EccHill elements is tested for other perturbations in addition to $J_2$.

\begin{table}
  \centering
  \caption{CPU times in milliseconds for high-order mapping (building the map, mapping and total using 5th-order Taylor expansions) and numerical propagation of a HEO orbit (test case 6) for 10,000 revolutions. $^\ast$Includes computing additional high-order map to solve for Poincar\'e map condition.}
    \begin{tabular}{lcccc}
    \hline
    Element set & \multicolumn{3}{c}{High-order mapping [ms]} & Numerical [ms] \\
          & Build map & Mapping & Total &  \\
    \hline
    COE   & 1061  & 16    & 1077  & 56784 \\
    MEE   & 2340  & 16    & 2356  & 51293 \\
    Hill  & 484   & 15    & 499   & 21294 \\
    EccHill & 624   & 16    & 640   & 33213 \\
    Cyl   & 2590$^\ast$  & 31    & 2621  & 20671 \\
    CylHz & 2278$^\ast$  & 31    & 2309  & 20982 \\
    Ideal & 3666  & 93    & 3759  & 25740 \\
    \hline
    \end{tabular}%
  \label{tab:ComputationTimesHEOJ2i30}%
\end{table}%

\subsection{Higher-order zonal and drag perturbations}
In the following test cases, we consider the perturbations due to $J_2$, $J_3$, $J_4$ and drag. These perturbations can be computed using the expressions in Sections \ref{sec:geopotential} and \ref{sec:drag}. In addition, the equations for calculating the velocity relative to the atmosphere in EccHill elements are given in Appendix~\ref{app:drag}. The LEO orbits considered here have an eccentricity of 0.01 and 0.001 and altitudes of 500 and 800 km, which are typical for Earth observation orbits. 

\subsubsection{Test case 7}
\begin{figure}
  \centering
	\includegraphics[width=0.8\columnwidth]{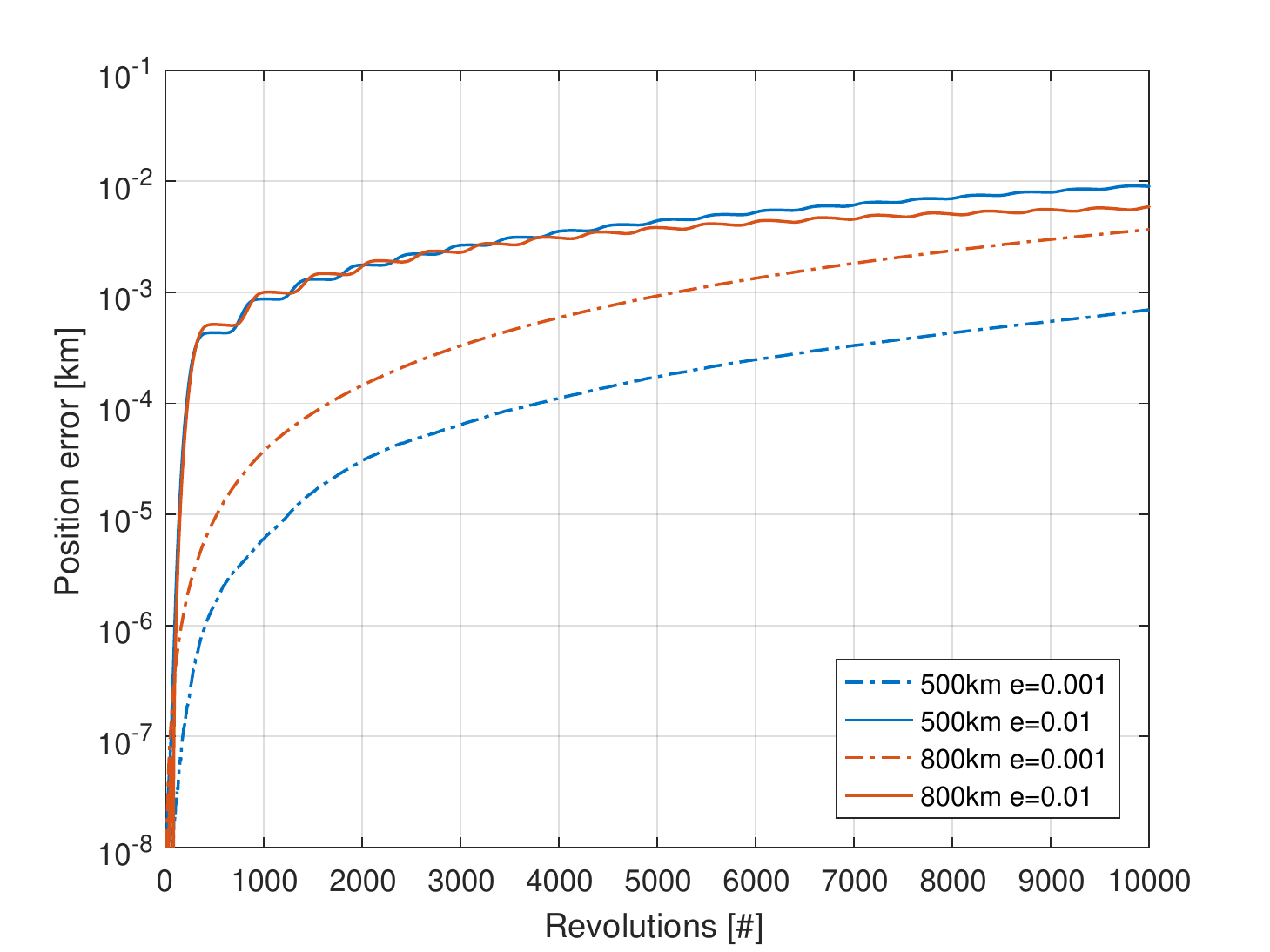}
    \caption{Position error using eccentric Hill variables for test case 7 (LEO, $i=30\degr$, $J_2-J_4$) for 10,000 revolutions at 500 km (655 days) and 800 km (698 days) altitude with $e=0.001$ and $e=0.01$.}
    \label{fig:LEO_J2J3J4_i30_10000rev_PosError}
\end{figure}

Figure~\ref{fig:LEO_J2J3J4_i30_10000rev_PosError} shows the position errors when considering zonal perturbations up to $J_4$ for 10,000 revolutions. The errors are similar to the $J_2$ only case (compare with Fig.~\ref{fig:500km_e01_i30_10000rev_PosError}) which means that the higher-order zonal perturbations do not affect the accuracy of the mapping. Besides, one can see that the mapping accuracy is higher for lower eccentricities, because the variation in the EccHill elements $\hat{f}$ and $\hat{g}$ is smaller for lower eccentricities.

\subsubsection{Test case 8}
\begin{figure}
  \centering
	\includegraphics[width=0.8\columnwidth]{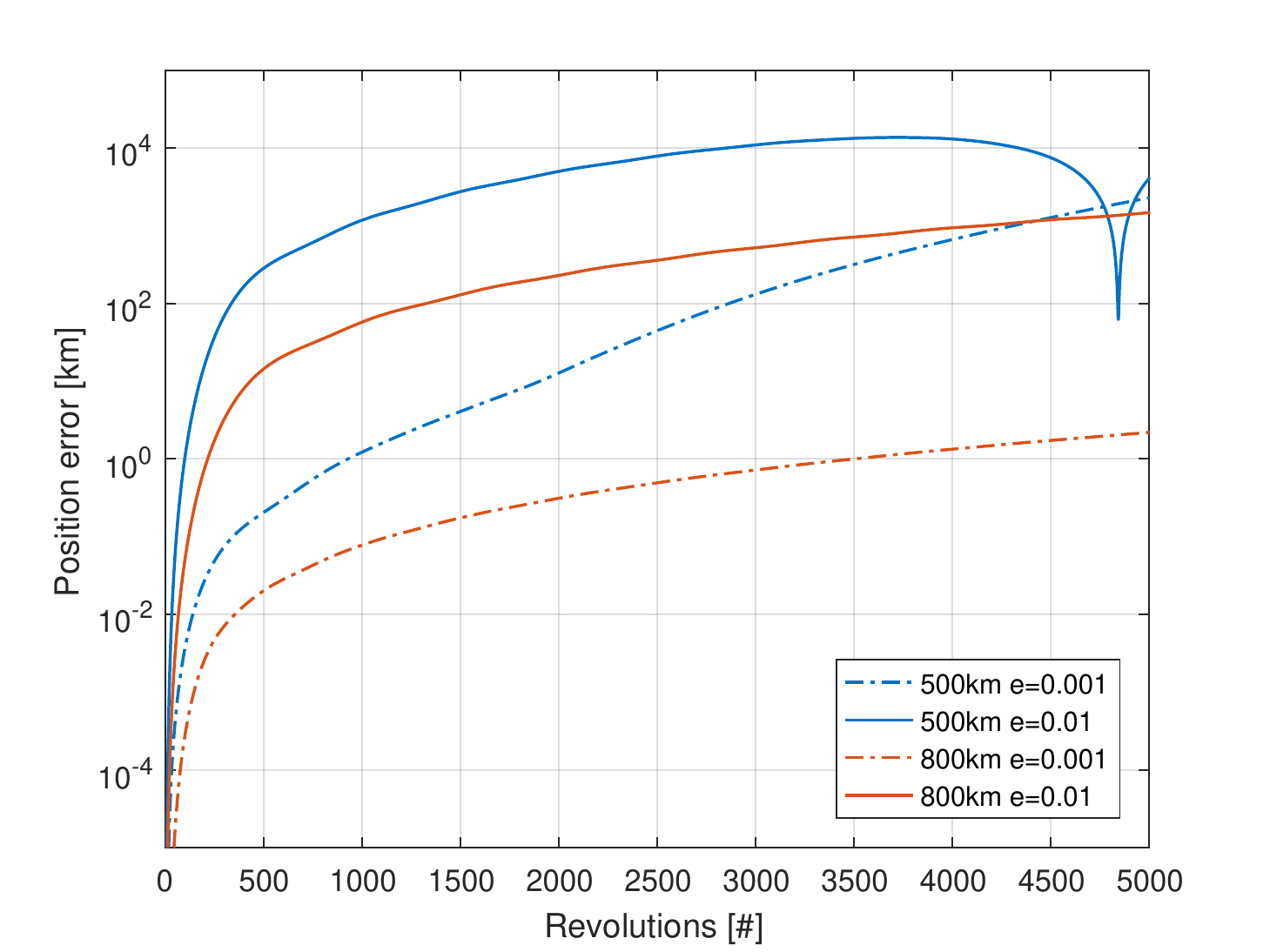}
    \caption{Position error using eccentric Hill variables for test case 8 (LEO, $i=30\degr$, $J_2-J_4$ and drag) for 5000 revolutions at 500 km (326 days) and 800 km (349 days) altitude with $e=0.001$ and $e=0.01$.}
    \label{fig:LEO_J2J3J4drag_i30_10000rev_PosError}
\end{figure}

\begin{figure}
  \centering
	\includegraphics[width=0.8\columnwidth]{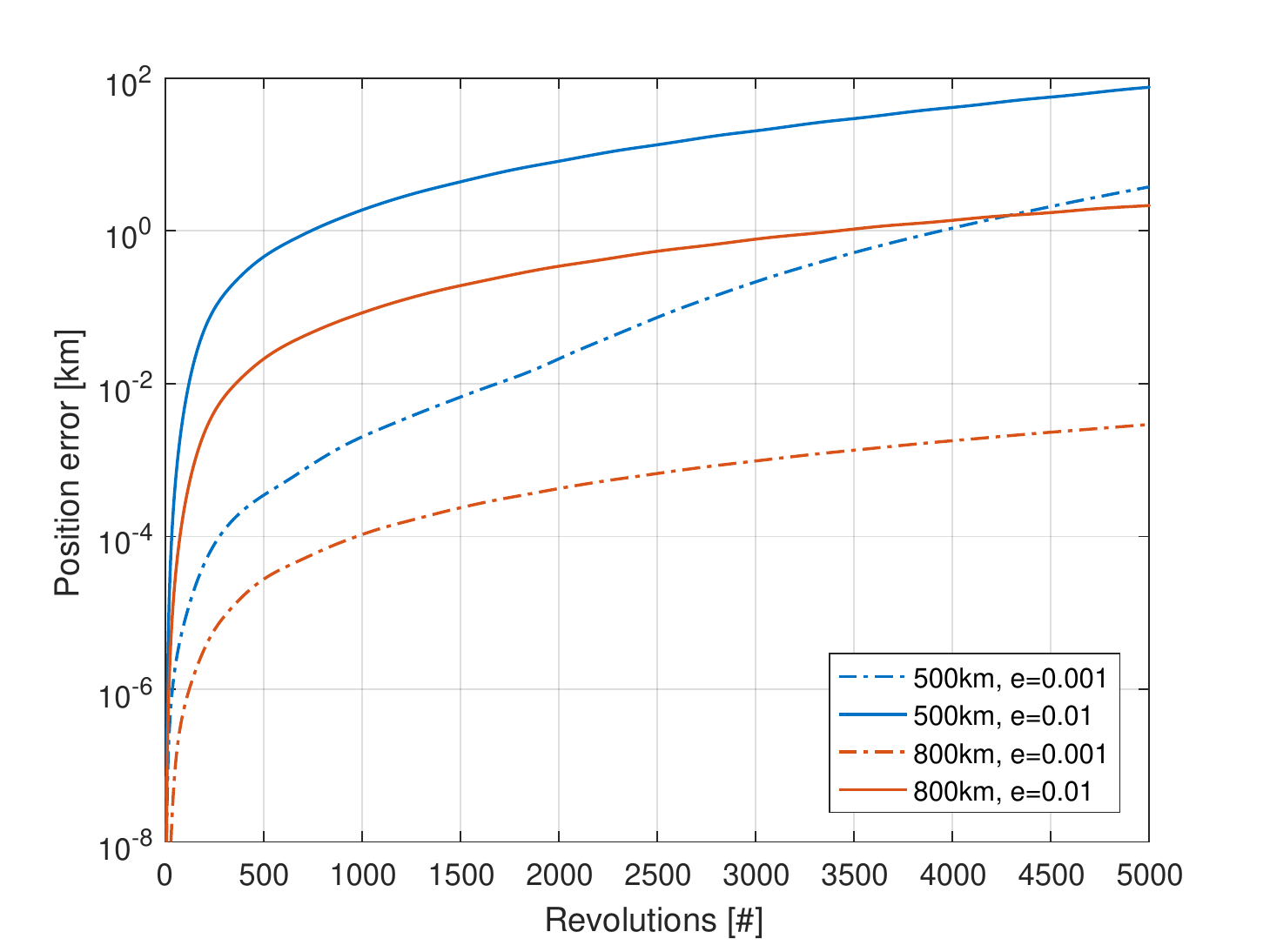}
    \caption{Position error without time error using eccentric Hill variables for test case 8 (LEO, $i=30\degr$, $J_2-J_4$ and drag) for 5000 revolutions at 500 km (326 days) and 800 km (349 days) altitude with $e=0.001$ and $e=0.01$.}
    \label{fig:LEO_J2J3J4drag_i30_10000rev_PosError_ElementsOnly}
\end{figure}

The accuracy of the mapping when drag is included is shown in Fig.~\ref{fig:LEO_J2J3J4drag_i30_10000rev_PosError}. The drag strongly affects the accuracy, because the angular momenta $H$ and $H_z$ now change secularly instead of being constant (on average). Nevertheless, for the orbits with $e=0.001$ at 500 and 800 km altitude the position error is less than 1 km for 60.6 days (926 revs) and 244.5 days (3501 revs), respectively. When the eccentricity is 0.01 the effect of drag is stronger and the mapping has an accuracy better than 1 km for 6.3 and 14.9 days at 500 and 800 km altitude, respectively. This means that the high-order map could be useful for fast on-board orbit prediction for several days or weeks, e.g. for navigation and control. Moreover, if the error in time is not considered, than the results including drag are much more accurate and the error remains below 1 km for extended periods of time especially when $e=0.001$, see Fig.~\ref{fig:LEO_J2J3J4drag_i30_10000rev_PosError_ElementsOnly}. 


\subsection{Fixed point}
This test case aims to illustrate how the HOTM technique can be used both to find a fixed point close to an initial guess and to study its center manifold. We use as a first guess a circular sun-synchronous orbit at 500 km altitude (orbital elements are provided Table~\ref{tab:testcases}) and we build a HOTM centred at this orbit in EccHill variables. Then we apply the approach introduced in Sec.\,\ref{sec:fixedPoint} to compute the nearby periodic orbit in the $J_2$-$J_4$ zonal problem. Starting with the initial guess $\hat{f}_0 = \hat{g}_0 = 0$, the numerical solver (\texttt{fsolve}) converges to the solution $\hat{f}_0 = 4.8222\times10^{-4}$ and $\hat{g}_0 = 1.0805\times10^{-3}$ (corresponding to $e=0.001183$ and $\omega=65.9489\degr$) in just four iterations. 

\begin{figure}
  \centering
	\includegraphics[width=0.8\columnwidth]{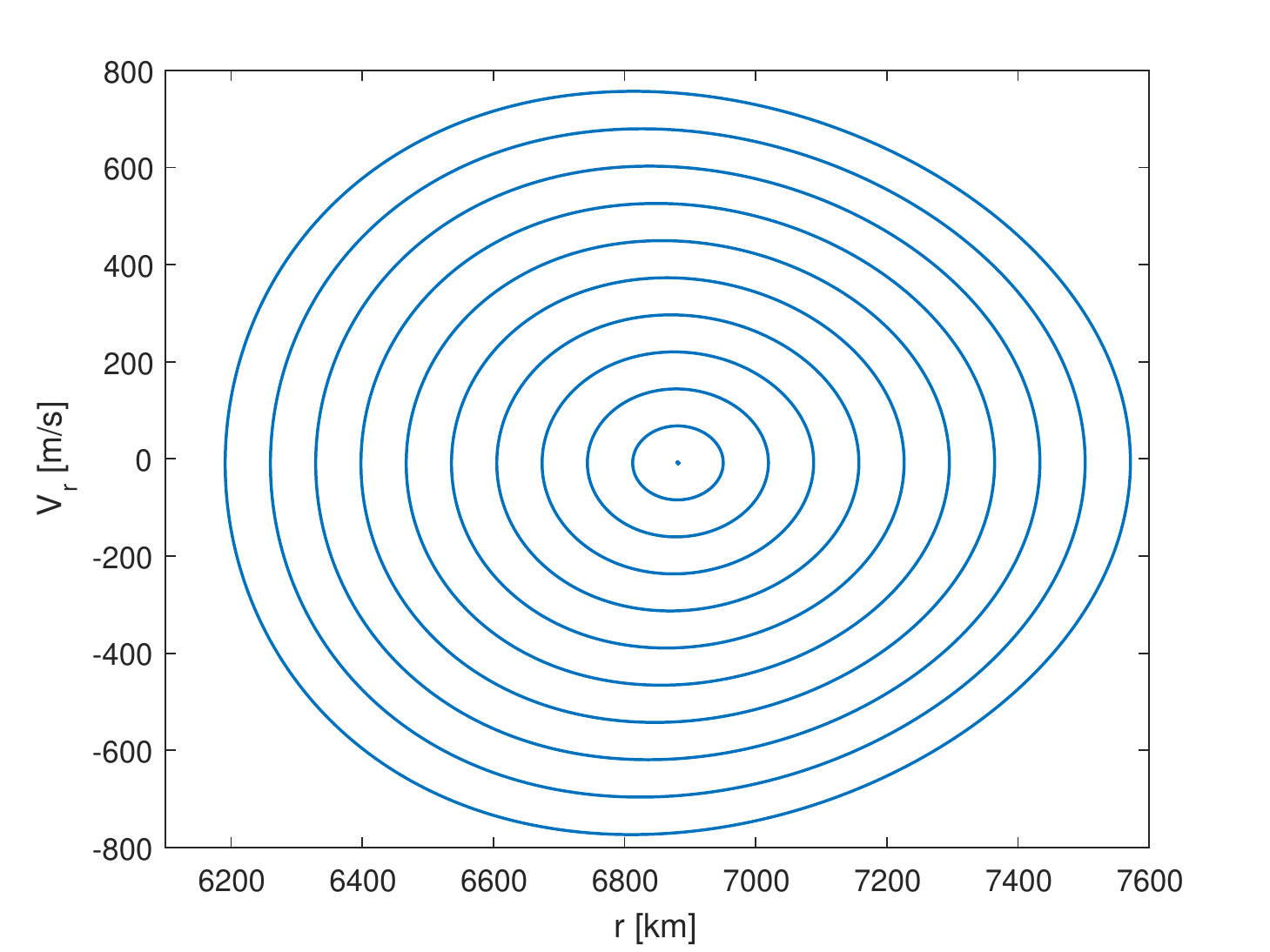}
    \caption{Poincar\'{e} section of radial distance and velocity $(r,V_r)$ at $z = 0$ for orbits with $E=-0.46365\mu/R_e$, $H_z=-0.13411\sqrt{\mu R_e}$ and $e=[0.001183,0.101183]$ perturbed by $J_2-J_4$.}
    \label{fig:FixedPointJ2J3J4_PoincareSection}
\end{figure}

After that, using the same high-order map we compute the invariant curves surrounding the periodic orbit parametrized in eccentricity, while keeping the orbital energy, $E$, and $H_z$ fixed. Figure~\ref{fig:FixedPointJ2J3J4_PoincareSection} shows the Poincar\'e section, in the $r,V_r$ space, for the periodic orbit (fixed point in the centre of the section) and quasi-periodic orbits (the invariant curves) with the eccentricity increased in steps of 0.01. The Poincar\'e plot was built using 2000 mappings, since the considered quasi-periodic orbits complete one cycle in less than 1560 revolutions. 

Noticeably, a single high-order map was used to compute the fixed point and construct the Poincar\'e surface of section. The maximum error in the mapped elements is $3.3\times10^{-7}$ km in position and 25.0 s in time and the combined error in position and time does not exceed 1 km up to an eccentricity of 0.041183. This demonstrates that using a single high-order map we can map the elements exactly and the time accurately for significant changes in eccentricity.


\begin{figure}
  \centering
	\includegraphics[width=0.8\columnwidth]{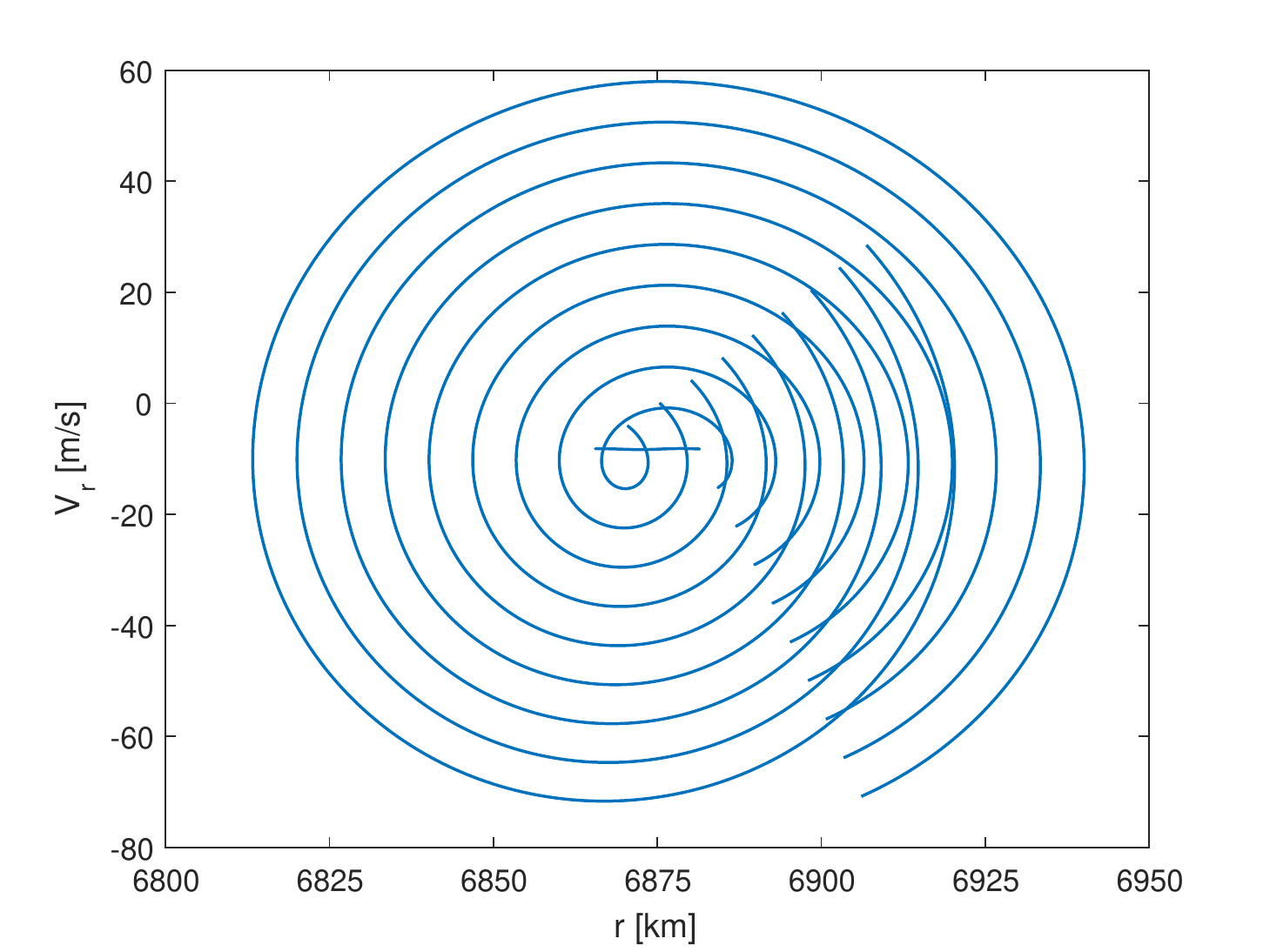}
    \caption{Poincar\'{e} section of radial distance and velocity $(r,V_r)$ at $z = 0$ for orbits with $E_0=-0.46365\mu/R_e$, $H_{z,0}=-0.13411\sqrt{\mu R_e}$ and $e_0=[0.001183,0.010183]$ perturbed by $J_2-J_4$ and drag for 2,000 revolutions (131 days).}
    \label{fig:FixedPointJ2J3J4drag_PoincareSection}
\end{figure}

Figure~\ref{fig:FixedPointJ2J3J4drag_PoincareSection} shows the Poincar\'e section when also drag is considered with an eccentricity up to 0.010183. One can see that the orbit contracts and circularises as the radial distance and velocity decrease due to drag. In addition, as one would expect, the ``fixed point'' is not fixed any more due to the decay in altitude. In this case, using a single high-order map to create the plot resulted in a maximum error of 0.73 km in position and 422 s in time after 2,000 mappings (131 days).

In summary, the proposed approach provides a tool for the efficient computation of Poincar\'e sections about stable fixed points. This enables qualitative dynamical system studies and has potential to be applied in challenging problems in astrodynamics such as the design and control of relative bounded motion, e.g. for formation flying missions.   

\subsection{Additional remarks}
\begin{figure}
  \centering
	\includegraphics[width=0.8\columnwidth]{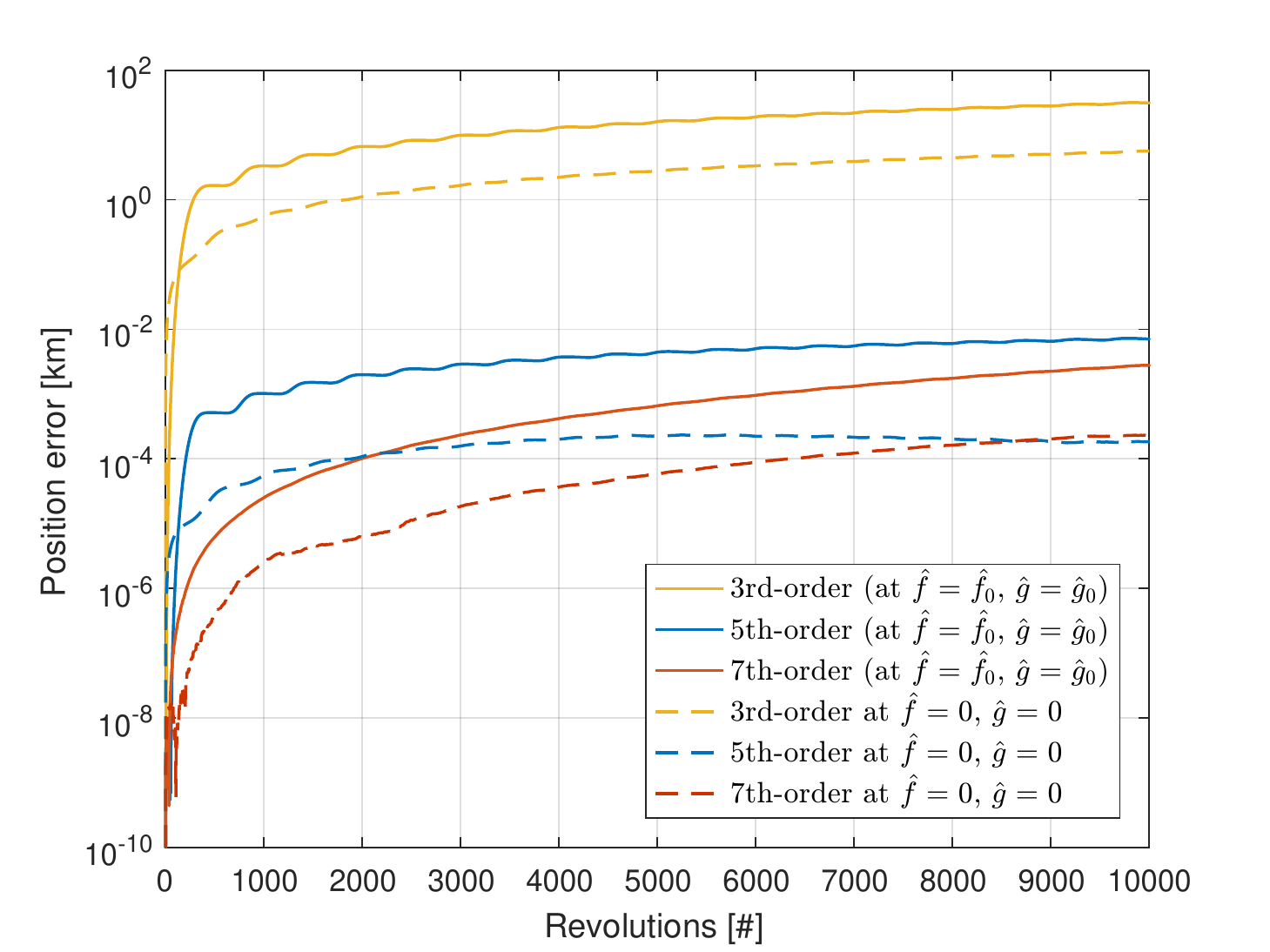}
    \caption{Position error using eccentric Hill variables with different expansion orders (3rd, 5th and 7th order) and by expanding $\hat{f}$ and $\hat{g}$ around their initial values (i.e. $\hat{f}_0$ and $\hat{g}_0$) or around zero for test case 1 (LEO, $i=30\degr$, $J_2$ only) for 10,000 revolutions  (655 days).}
    \label{fig:MidWhit_3rd5th7thO_expandAroundFGzero}
\end{figure}

The results shown in the paper were computed using 5th-order Taylor expansions. This expansion order provides a good trade-off between accuracy of the results and speed of the calculations. However, higher or lower expansion orders may be selected for improved accuracy or efficiency. Figure~\ref{fig:MidWhit_3rd5th7thO_expandAroundFGzero} shows the position error for test case 1 with EccHill elements using different expansion orders, namely 3rd, 5th and 7th order. As expected, the error reduces with increasing expansion order. On the other hand, the computation time increases from 63 ms to 281 ms and 1092 ms for 3rd, 5th and 7th order, respectively.

Moreover, the accuracy of the results can be improved by cleverly choosing the expansion point. For example, for the eccentric Hill variables $\hat{f}$ and $\hat{g}$ the variation is known to be in the domain $[-e,e]$. Therefore, the maximum deviation from the expansion point is smallest when the expansion point is centered in the domain, i.e. at zero. Figure~\ref{fig:MidWhit_3rd5th7thO_expandAroundFGzero} shows the position error when expanding $\hat{f}$ and $\hat{g}$ around zero using different expansion orders. The first tens of revolutions the position error is larger compared to expanding around initial values $\hat{f}_0$ and $\hat{g}_0$ because the deviation from the expansion point is larger. However, on the long term the position error is about one order magnitude smaller. This improvement in accuracy comes without any cost in computation time.

\section{Conclusions}
The performance of Poincar\'e and stroboscopic mapping using the high-order transfer map method was tested using different element sets for various types of orbits. The choice of coordinates has a strong effect on the accuracy of the mapped states. By choosing the proper elements, accurate mapping of the osculating orbital state under $J_2$ perturbation is possible for thousands of revolutions. Moreover, high-order mapping is an order of magnitude faster than numerically propagating an orbit for 10,000 revolutions and three orders of magnitude more efficient if the map is precomputed. The main causes of error growth are the drift in $\omega$ and $r$, and variation in $\Omega$ when coupled with the inclination.

A new set of orbital elements, the eccentric Hill variables, was introduced for improved high-order mapping accuracy. The elements can be used straightforwardly for both stroboscopic mapping in general and Poincar\'e mapping on the equatorial plane. Using the new elements, highly-elliptical orbits and orbits perturbed by higher-order zonal perturbations ($J_2$-$J_4$) can be mapped extremely accurately with only an error in the mapped time. Furthermore, the high-order mapping was shown to be accurate for drag-perturbed orbits for several weeks or months depending on the strength of the drag. 

As an example application, we used the method to compute a fixed point under $J_2$-$J_4$ perturbations and investigate the quasi-periodic orbits in a large domain around the fixed point (including drag). This was achieved using a single high-order map, which shows the potential of the method.

The high-order mapping technique could be applied to various problems in celestial mechanics that require the mapping of orbits, e.g. investigating the stability of orbits using Poincar\'e maps. Furthermore, in the future work the presented method can be extended to other perturbations, such as third-body and tesseral perturbations.

\begin{acknowledgements}
The authors are very grateful for the help, support and advice from Dr. Alexander Wittig on this work. 
David Gondelach is funded by the Surrey Space Centre (SSC) to undertake his PhD at the University of Surrey. Finally, the use of the Differential Algebra Computational Engine (DACE) developed by Dinamica Srl is acknowledged.
\end{acknowledgements}

\noindent \small \textbf{Conflicts of interest} The authors declare that there are no conflicts of interest regarding the publication of this paper.


\appendix

\section{Equations of motion}

\subsection{Gauss' equations of motion in modified equinoctial elements}
\label{app:EoMMEE}
The Gaussian equations of motion in terms of modified equinoctial elements are given by \citep{Walker1985,Walker1986}:
\begin{align}
	\frac{dp}{dt} =& \frac{2p}{w}\sqrt{\frac{p}{\mu}}f_{\theta}  ,	\\
	\frac{df}{dt} =& \sqrt{\frac{p}{\mu}} \left\{	f_r\sin{L}+\left[(w+1)\cos{L}+f\right]\frac{f_{\theta}}{w} \right. \\
	&~\left. -~(h\sin{L}-k\cos{L})\frac{g}{w}f_n \right\} ,	\\
	\frac{dg}{dt} =& \sqrt{\frac{p}{\mu}} \left\{ -f_r\cos{L}+\left[(w+1)\sin{L}+g\right]\frac{f_{\theta}}{w} \right. \\
	&~\left. +~(h\sin{L}-k\cos{L})\frac{f}{w}f_n \right\} ,	\\
	\frac{dh}{dt} =& \sqrt{\frac{p}{\mu}} \frac{s^2f_n}{2w}\cos{L} ,	\\
	\frac{dk}{dt} =& \sqrt{\frac{p}{\mu}} \frac{s^2f_n}{2w}\sin{L} ,	\\
	\frac{dL}{dt} =& \sqrt{\mu p}\left( \frac{w}{p} \right)^2 + \frac{1}{w}\sqrt{\frac{p}{\mu}} (h\sin{L}-k\cos{L})f_n	.
\end{align}\\

\subsection{Gauss' equations of motion in Hill variables}
\label{app:EoMpolarnodal}
Gauss' planetary equations can be written in Hill variables as follows \citep{mazzini2015flexible}:
\begin{align}
	\frac{dr}{dt} &= \dot{r}  ,	\\
	\frac{du}{dt} &= \frac{H}{r^2} - \frac{r \cos{i} \sin{u}}{H \sin{i}} f_n ,	\\
	\frac{d\Omega}{dt} &= \frac{r\sin{u}}{H \sin{i}}f_n ,	\\
	\frac{d\dot{r}}{dt} &= -\frac{\mu}{r^2} + \frac{H^2}{r^3} + f_r ,	\\
	\frac{dH}{dt} &= r f_t ,	\\
	\frac{dH_z}{dt} &= r \cos{i} f_t - r\sin{i}\cos{u}  f_n	.
\end{align}\\

\section{Velocity in eccentric Hill variables}
\label{app:drag}
The components of the velocity vector of the satellite $\mathbfit{V}$ in the orbital frame (see \citet{battin1999introduction}) can be expressed in eccentric Hill variables as follows:
\begin{align}
	V_r &= \frac{\mu}{H} e\sin{\nu} = \frac{\mu}{H} (\hat{f}\sin{u} - \hat{g}\cos{u}) ,
	\notag \\
	V_t &= \frac{\mu}{H} (1+e\cos{\nu}) = \frac{\mu}{H} (1+\hat{f}\cos{u} + \hat{g}\sin{u}) ,
	\notag \\
	V_h &= 0 .
\end{align}
The atmosphere is assumed to rotate with the Earth about its $z$-axis with an angular velocity $\omega_e$.
The velocity with respect to the atmosphere $\mathbfit{V}_{rel}$ in the orbital frame is then obtained by subtracting the local velocity of atmosphere \citep{vallado2013fundamentals}:
\begin{equation}
	\mathbfit{V}_{rel} = \mathbfit{V} - \boldsymbol{\omega}_e \times \mathbfit{r} = \begin{bmatrix}
	\frac{\mu}{H} (\hat{f}\sin{u} - \hat{g}\cos{u}) \\
	\frac{\mu}{H} (1+\hat{f}\cos{u} + \hat{g}\sin{u}) - r\,\omega_e\cos{i} \\
	r\,\omega_e \sin{i} \cos{u}
	\end{bmatrix}
\end{equation}


\bibliographystyle{spbasic}      
\bibliography{library}   

\end{document}